# Graph-based Logic and Sketches I: The General Framework

Atish Bagchi[*] and Charles Wells

October 12, 1994

## 1 Introduction

Sketches as a method of specification of mathematical structures are an alternative to the string-based specification employed in mathematical logic. Sketch theory has been criticized as being lacunary when contrasted with logic because it apparently has nothing corresponding to proof theory. This article describes an approach to associating a proof-theoretic structure with a sketch. It is written in an innovative style using some of the ideas for presenting mathematics outlined in [Bagchi and Wells, 1993b].

### 1.1 Prerequisites

This article requires familiarity with the basic notions of mathematical logic as in Chapters 2 through 5 of [Ebbinghaus *et al.*, 1984], and with category theory and sketches as in Chapters 1 through 10 of [Barr and Wells, 1990]. We specifically presuppose that finite-limit sketches are known. Some notation for these ideas is established in Section 4.

[*]The first author is grateful to the Mathematical Sciences Research Institute for providing the opportunity to work on this paper there. Research at MSRI is supported in part by NSF grant #DMS-9022140.

## 1.2 Rationale

Traditional treatments of formal logic provide:

1. A syntax for formulas. The formulas are typically defined recursively by a production system (grammar).

2. Inference relations between sets of formulas. This may be given by structural induction on the construction of the formulas.

3. Rules for assigning meaning to formulas (semantics) that are sound with respect to the inference relation. The semantics may also be given by structural induction on the construction of the formulas.

First order logic, the logic and semantics of programming languages, and the languages that have been formulated for various kinds of categories are all commonly described in this way. The formulas in those logics are strings of symbols that are ultimately modeled on the sentences mathematicians speak and write when proving theorems.

Mathematicians with a category-theoretic point of view frequently state and prove theorems using graphs and diagrams. The graphs, diagrams, cones and other data of a sketch are formal objects that correspond to the graphs and diagrams used by such mathematicians in much the same way as the formulas of traditional logic correspond to the sentences mathematicians use in proofs. The functorial semantics of sketches is sound in the informal sense that it preserves by definition the structures given in the sketch. The analogy to the semantics of traditional model theory is close enough that sketches and their models fit the definition of "institution" ([Goguen and Burstall, 1986]), which is an abstract notion of a logical system having syntactic and semantic components. This is described in detail in [Barr and Wells, 1990], Section 10.3. Note that the soundness of functorial semantics appears trivial when contrasted with the inductive proofs of soundness that occur in string-based logic because the semantics functor is not defined recursively.

## 1.3 Overview

This paper exhibits a structure in sketch theory that corresponds to items 1 and 2 in the description of logic in Section 1.2. The data making up the structure we give do not correspond in any simple way to the data involved in items 1 and 2 of traditional logic; we will discuss the differences *in situ*.

The data in this structure are parametrized by the type of categorial theory being considered. Let $E$ denote a finite-limit sketch that presents a type of category as essentially algebraic over the theory of categories. Let **E** denote the finite-limit theory generated by $E$. **E**-categories are then the models of **E** in the category of sets. (This is described in Section 5.) The kinds of categories that can be described in this way include categories with finite products, categories with limits or colimits over any particular set of diagrams, cartesian closed categories, regular categories, toposes, and many others. Sketches for several specific instances of **E** are given in Appendix B.

An **E**-sketch $S$ (defined in Section 5.3.1 below) is a sketch that allows the specification of any kind of construction that can be made in an **E**-category. As an example, let **CCC** be a finite-limit theory for Cartesian closed categories. It is possible to require that a certain object in a **CCC**-sketch be the function space $A^B$ of two other objects $A$ and $B$ of the sketch. This generalization greatly enhances the expressive power of sketches as originally defined by Ehresmann (in which only limits and colimits can be specified).

Such an **E**-sketch is described in Section 5.3.1 as a formally adjoined global element of the limit vertex of a certain diagram in **E**. This diagram forces the value of the limit vertex in a model to contain the graph, diagrams and other structures in the sketch. The **E**-theory of this sketch, which corresponds roughly, but not precisely, to the deductive closure of a list of axioms in traditional logic, is the initial model of **E**[$S$], which is **E** with the global element adjoined.

An assertion in this setting is a potential factorization (PF) of an arrow of the **E**-theory of $S$ through an arrow into its codomain (defined precisely in Section 6.1). The assertion is valid if the PF does indeed factorize in every model of **E**[$S$].

Instead of the set of rules of deduction of a traditional theory, we have a set of rules of construction. More precisely, we give in Appendix A a system of construction rules that produce the categorial theory of a finite-limits sketch. We say that the potential factorization is **deducible** if there is an actual factorization in the categorial theory of the sketch. Such an arrow must be constructible by the rules in Appendix A. Thus the usual system of inference is replaced by a system of construction of arrows in a finite-limit category (no matter what type of category is sketched by $E$). This system is sound and complete with respect to models (Section 6.2).

The fact that we have assumed finite-limit sketches as given prior to the general definition of **E**-sketch is basic to the strategy of the paper, which is to make finite-limit logic the basic logic for all **E**-sketches (Section 8). The variation in what can be proved, for example for finite-limit sketches (**E**-sketches where **E** = **FinLim**) as contrasted with cartesian closed sketches (in which case **E** is a finite-limit theory for cartesian closed categories such as **CCC** in Appendix B.6) is entirely expressed by the choice of **E** and has no effect on the rules of construction.

The system that we describe may not facilitate human manipulation. However, as we discuss in Remark A.1.2, the system should not be difficult for computer implementation. Our approach may seem unusual from the point of view of symbolic logic, but it is in keeping with practices in computer implementation of complex data structures (a proof is a complex data structure) in that much more detail about the relationships between different parts of the structure are necessarily made explicit than is customary for systems that are designed to be understood by human beings, who have sophisticated built-in pattern recognition abilities.

### 1.4 Acknowledgments

This article is better because of conversations we have had with Robin Cockett and Colin McLarty. The names "string-based logic" and "graph-based logic" were suggested by Peter Freyd. The diagrams were prepared using K. Rose's `xypic`.

## 2 Notation

The notation established in this section is used throughout the article.

### 2.1 Lists

Given a set $A$, $\mathsf{List}[A]$ denotes the set of lists of elements of $A$, including the empty list. The $k$th entry in a list $w$ of elements of $A$ is denoted by $w_k$ and the length of $w$ is denoted by $\mathsf{Length}[w]$. The **range** of $w$, denoted by $\mathsf{Rng}[w]$, is the set of elements of $A$ occurring as entries in $w$. If $f : A \to B$ is a function, $\mathsf{List}[f] : \mathsf{List}[A] \to \mathsf{List}[B]$ is by definition $f$ "mapped over" $\mathsf{List}[A]$: If $w$ is a list of elements of $A$, then the $k$th entry of $\mathsf{List}[f](w)$ is by definition $f(w_k)$. This makes $\mathsf{List}$ a functor from the category of sets to itself.

## 2.2 Graphs

For a graph $G$, the set of nodes of $G$ is denoted by $\mathsf{Nodes}[G]$ and the set of arrows is $\mathsf{Arrows}[G]$. The underlying graph of a category $\mathscr{C}$ is denoted by $\mathsf{UndGr}[\mathscr{C}]$. A subgraph $H$ of a graph $G$ is said to be **full** if every arrow $f : h_1 \to h_2$ of $G$ between nodes of $H$ is an arrow of $H$.

## 2.3 Diagrams

We define an equivalence relation on the set of graph homomorphisms into a graph $G$ as follows: If $\delta : I \to G$ and $\delta' : I' \to G$ are graph homomorphisms, then $\delta$ is equivalent to $\delta'$ if and only if there is a graph isomorphism $\phi : I \to I'$ such that

$$\begin{array}{ccc} I & \xrightarrow{\phi} & I' \\ & \searrow{\delta} \quad \swarrow{\delta'} & \\ & G & \end{array} \qquad (1)$$

commutes. A **diagram** in $G$ is then by definition an equivalence class of graph homomorphisms $\delta : I \to G$. As is the practice when an object is defined to be an equivalence class, we will refer to a diagram by any member of the equivalence class.

If $\delta : I \to G$ is a diagram, $I$ is said to be a **shape graph** of the diagram, denoted by $\mathsf{ShpGr}(\delta)$, and $G$ is said to be the **ambient space** of the diagram, denoted by $\mathsf{AmbSp}(\delta)$. Observe that the ambient space of the diagram is determined absolutely, but the shape graph is determined only up to an isomorphism that makes Diagram (1) commute.

If $I$ is a graph, $\mathscr{C}$ is a category, and $\delta : I \to \mathsf{UndGr}[\mathscr{C}]$ is a diagram, we will also write it as $\delta : I \to \mathscr{C}$ (unless there is danger of confusion).

### 2.3.1 Convention on drawing diagrams

It is customary to draw a diagram without naming its shape graph. We adopt the following convention: If a diagram is represented by a drawing on the page, the shape graph of the diagram is the graph that has one node for each object shown and one arrow for each arrow shown, with source and target as shown. Two objects at different locations on the page correspond to two different nodes of the shape graph, *even if they have the same label*, and an analogous remark applies to arrows.

### 2.3.2 Example

The diagram (2)

$$A \xrightarrow{f} A$$
$$g \downarrow \quad \downarrow g \qquad (2)$$
$$B \xrightarrow{\mathsf{Id}[B]} B$$

called $\delta$, has shape graph

$$h \xrightarrow{t} i$$
$$u \downarrow \quad \downarrow v \qquad (3)$$
$$j \xrightarrow{x} k$$

so that $\delta(h) = A$, $\delta(i) = A$, $\delta(v) = g$, $\delta(x) = \mathsf{Id}[B]$ and so on. Diagram (4) below also has shape graph (3) (or one isomorphic to it, of course):

$$A \xrightarrow{f} B$$
$$g \downarrow \quad \downarrow h \qquad (4)$$
$$C \xrightarrow{k} D$$

On the other hand, Diagram (5) below

$$A \xrightarrow{f} A$$
$$g \searrow \quad \swarrow g \qquad (5)$$
$$B$$

*is not the same diagram as (2).* It has shape graph

$$i \xrightarrow{u} j$$
$$v \searrow \quad \swarrow w \qquad (6)$$
$$k$$

## 2.4 Cones

For any graph $G$ and diagram $\delta : I \to G$, a cone $\Theta : v \lhd (\delta : I \to G)$ (also written $\Theta : v \lhd \delta$ if the context makes this clear) has **vertex** $v$ denoted by Vertex[Θ] and **base diagram** $\delta$ denoted by BsDiag[Θ]. For each node $i$ of ShpGr[$\delta$], the formal projection of the cone Θ from Vertex[Θ] to $\delta(i)$ is denoted by Proj[Θ, $i$]: $v \to \delta(i)$. For a category $\mathscr{C}$, a cone $\Theta : v \lhd (\delta : I \to \mathscr{C})$ is **commutative** if, for every arrow $f : i \to j$ in $I$, the diagram

$$\begin{array}{c} \text{Vertex}[\Theta] \\ \text{Proj}[\Theta,i] \swarrow \quad \searrow \text{Proj}[\Theta,j] \\ \delta(i) \xrightarrow{\delta(f)} \delta(j) \end{array} \tag{7}$$

commutes. The limit cone for a diagram $\delta$ will be denoted by LimCone[$\delta$] : Lim[$\delta$] $\lhd$ $\delta$.

### 2.5 Fonts

In general, variable objects are given in slant or script notation and specific objects (given by proper names) are given in upright notation. In more detail, we have the following systematic notation.

1. Specific data constructors, such as List, and specific fieldnames for complex objects, such as Nodes[$G$], are given in sans serif and are capitalized as shown.

2. Specific objects and arrows of sketches are also given in sans serif.

3. Specific constructor spaces, such as **FinLim** and **CCC**, are given in **bold sans serif**. We use **E** to denote a variable constructor space because of the unavailability of bold slanted sans serif.

4. Specific categories other than constructor spaces, such as **Set**, are given in **boldface**.

5. Diagrams (specific and variable) are named by lowercase Greek letters.

6. Cones (specific and variable) are named by uppercase Greek letters.

7. Models (specific and variable) are given in uppercase fraktur, for example $\mathfrak{M}$, $\mathfrak{C}$.

8. Variable sketches, and objects and arrows therein, are given in *slanted sans serif*.

9. Variable categories other than constructor spaces are given in script, for example $\mathscr{A}$, $\mathscr{B}$, $\mathscr{C}$.

10. Other variable objects are given in math italics, such as $a$, $b$, $c$ or (especially arrows) in lowercase Greek letters.

# 3 Adjoining nodes to diagrams

We now state and prove some technical lemmas that are used many times herein.

## 3.1 Restrictions of diagrams

**3.1.1 Definition** Let $\delta: I \to \mathscr{C}$ be a diagram and $\mathsf{Incl}[J \subseteq I]: J \to I$ an inclusion of graphs. Let $\Theta: v \lhd \delta$ be a cone. The **base-restriction of $\Theta$ to $J$** is defined to be the cone $\Theta|^J : v \lhd (\delta \circ \mathsf{Incl}[J \subseteq I])$ with vertex $v$ and projections defined by $\mathsf{Proj}[\Theta|^J, j] := \mathsf{Proj}[\Theta, j]: v \to \delta(j)$ for all nodes $j$ of $J$. In this case, we also say that $\Theta$ is a **base-extension** of $\Theta|^J$.

**3.1.2 Remark** Observe that if $\Theta$ is commutative, then so is $\Theta|^J$.

**3.1.3 Definition** Let $\delta: I \to \mathscr{C}$ be a diagram and $J$ a subgraph of $I$. Then the subdiagram $\delta|_J$ is said to **dominate** $\delta$, or to be **dominant in $\delta$**, if every commutative cone $\Theta: v \to (\delta|_J)$ has a base extension to a commutative cone $\Theta': v \to \delta$ with the same vertex.

**3.1.4 Remark** Tholen and Tozzi [1989] give a condition ("confinality") on $I$ and $J$ such that any diagram based on $I$ is dominated by its base restriction to $J$. One type of dominance that their condition does not cover is the case in which $\delta$ is obtained from $\delta|_J$ by adjoining a limit cone over a subdiagram of $\delta|_J$. This occurs in Section 7.1.1.

## 3.2 Limits of subdiagrams

**3.2.1 Lemma** *Let $\delta: I \to \mathscr{C}$ be a diagram and let $J$ be a subgraph of $I$ with inclusion $\mathsf{Incl}[J \subseteq I]$. Let $\gamma = \delta|_J$. Then there is a unique arrow*

$\phi : \mathsf{Lim}[\delta] \to \mathsf{Lim}[\gamma]$ *such that for all nodes $j$ of $J$,*

$$
\begin{array}{c}
\mathsf{Lim}[\delta] \xrightarrow{\phi} \mathsf{Lim}[\gamma] \\
\mathsf{Proj}[\mathsf{LimCone}[\delta], j] \searrow \swarrow \mathsf{Proj}[\mathsf{LimCone}[\gamma], j] \\
\delta(j)
\end{array}
\quad (8)
$$

*commutes. Moreover, if $\gamma$ is a dominant subdiagram of $\delta$, then $\phi$ is an isomorphism.*

**Proof** The base-restriction of $\mathsf{LimCone}[\delta]$ to $J$ is a commutative cone $\Theta$ with vertex $\mathsf{Lim}[\delta]$. We choose $\phi$ to be $\mathsf{Fillin}[\Theta, \gamma]$, the arrow produced by Rule ∃FIA of Appendix A. It follows from Rule CFIA that for each node $j$ of $J$, the diagram below

$$
\begin{array}{c}
\mathsf{Lim}[\delta] \xrightarrow{\mathsf{Fillin}[\Theta,\gamma]} \mathsf{Lim}[\gamma] \\
\mathsf{Proj}[\Theta, j] \searrow \swarrow \mathsf{Proj}[\mathsf{LimCone}[\gamma], j]] \\
\gamma(j)
\end{array}
\quad (9)
$$

commutes. Because they do, rule !FIA implies that $\phi$ is the only arrow that makes them all commute, as required.

Assume that $\gamma$ dominates $\delta$. Then we may choose an extension $\Psi : \mathsf{Lim}[\gamma] \to \delta$ of $\mathsf{LimCone}[\gamma]$ to $\delta$. Using ∃FIA, we define

$$\psi := \mathsf{Fillin}[\Psi, \delta] : \mathsf{Lim}[\gamma] \to \mathsf{Lim}[\delta]$$

It follows from !FIA that $\psi$ is the only arrow from $\mathsf{Lim}[\gamma]$ to $\mathsf{Lim}[\delta]$ that makes all diagrams of the form

$$
\begin{array}{c}
\mathsf{Lim}[\gamma] \xrightarrow{\psi} \mathsf{Lim}[\delta] \\
\mathsf{Proj}[\Psi, j] \searrow \swarrow \mathsf{Proj}[\mathsf{LimCone}[\delta], j]] \\
\gamma(j)
\end{array}
\quad (10)
$$

commute for each node $j$ of $J$. Since $\phi \circ \psi : \mathsf{Lim}[\gamma] \to \mathsf{Lim}[\gamma]$ and $\mathsf{Id}[\mathsf{Lim}[\gamma]] : \mathsf{Lim}[\gamma] \to \mathsf{Lim}[\gamma]$ both commute with all the projections, it follows from !FIA that $\phi \circ \psi = \mathsf{Id}[\mathsf{Lim}[\gamma]]$, A similar argument shows that $\psi \circ \phi = \mathsf{Id}[\mathsf{Lim}[\delta]]$, so that $\phi$ is an isomorphism. □

## 3.3 Special cases of extending diagrams

Here we define some special cases of dominance that are easy to recognize.

**3.3.1 Definition** Let $I$ and $J$ be graphs that have the same nodes, and suppose that $I$ has exactly one arrow $a : j \to k$ not in $J$. Let $\delta : I \to \mathscr{C}$ be a diagram with the property that for all nodes $j'$ of $J$ and all arrows $f : j \to j'$ and $g : j' \to k$,

$$\begin{array}{ccc} j & \xrightarrow{a} & k \\ & \searrow_{f} \nearrow_{g} & \\ & j' & \end{array} \tag{11}$$

commutes. Then we say $\delta$ **extends** $\delta|_J$ **by adjoining a composite**.

**3.3.2 Definition** Let $I$ and $J$ be graphs with the following properties:

AC.1 $I$ has exactly one node $v$ not in $J$.

AC.2 Every arrow in $I$ not in $J$ has target $v$.

Suppose that $\delta : I \to \mathscr{C}$ is a diagram with the property that if $a : i \to v$, $b : j \to v$ and $f : i \to j$ are arrows of $I$, then

$$\begin{array}{ccc} i & \xrightarrow{a} & v \\ & \searrow_{f} \nearrow_{b} & \\ & j & \end{array} \tag{12}$$

commutes. Then we say $\delta$ **extends** $\delta|J$ by **adjoining a commutative cocone**.

**3.3.3 Definition** Let $I$, $J$ and $J'$ be graphs with $J' \subseteq J \subseteq I$, such that $J'$ is full in $J$, $I$ contains exactly one node $v$ not in $J$, and for each node $j$ of $J'$, $I$ contains exactly one arrow $p_j : v \to j$ and no other arrows not in $J$. Let $\delta : J \to \mathscr{C}$ be a diagram, and suppose $\delta' : I \to \mathscr{C}$ is a diagram that extends $\delta$ in such a way that $\delta'(v)$ and the arrows $\delta'(p_j)$ constitute a limit cone to $\delta|_{J'}$. Then we say that $\delta'$ **extends** $\delta$ **by adjoining a limit**.

**3.3.4 Lemma**  *Suppose that $\delta' : I \to \mathscr{C}$ extends $\delta : J \to \mathscr{C}$ by adjoining a composite, a commutative cocone, or a limit. Then*

$$\mathsf{Fillin}[\mathsf{LimCone}[\delta']|_J, \delta] : \mathsf{Lim}[\delta'] \to \mathsf{Lim}[\delta]$$

*is an isomorphism.*

**Proof**  It is easy to see that, in all three cases, $J$ is dominant in $I$. □

## 4 Finite limit sketches

We use a general concept of sketch described in Section 5 that is based on the concept of finite limit sketch, a particular case of projective sketch due to Ehresmann. In this section, we review briefly some aspects of finite limit sketches that are relevant later.

### 4.1 Categorial theories

The categorial theory generated by a finite-limit sketch $S$ is denoted by $\mathsf{CatTh}[\mathbf{FinLim}, S]$. It is a category with finite limits together with a sketch morphism

$$\mathsf{UnivMod}[\mathbf{FinLim}, S] : S \to \mathsf{CatTh}[\mathbf{FinLim}, S]$$

called the **universal model** of the sketch. It has the following property: For each model $\mathfrak{M}$ of $S$, there exists a functor $F_\mathfrak{M} : \mathsf{CatTh}[\mathbf{FinLim}, S] \to \mathbf{Set}$, determined uniquely up to natural isomorphism, with the property that

$$U(F_\mathfrak{M}) \circ \mathsf{UnivMod}[\mathbf{FinLim}, S] = \mathfrak{M}$$

where $U$ is the forgetful functor from the category of finite-limit categories and finite-limit-preserving functors to the category of finite-limit sketches and morphisms thereof. It follows from these properties that $\mathsf{CatTh}[\mathbf{FinLim}, S]$ is determined up to equivalence of categories and $\mathsf{UnivMod}[\mathbf{FinLim}, S]$ is determined up to natural isomorphism.

In this article, for a given sketch $S$, we assume given a particular instance of $\mathsf{CatTh}[\mathbf{FinLim}, S]$: that constructed in [Barr and Wells, 1993a]. It has the following properties (which are not preserved by equivalence of categories):

T.1  $\mathsf{CatTh}[\mathbf{FinLim}, S]$ is a category with canonical finite limits. (The construction explicitly produces the canonical limits.)

T.2 Every object of CatTh[**FinLim**, S] is the canonical limit of a diagram of the form UnivMod[**FinLim**, S]∘δ where δ is a diagram in the graph of S.

T.3 Every arrow of CatTh[**FinLim**, S] is a composite of projections from canonical limits and arrows of the form UnivMod[**FinLim**, S]($f$) for arrows $f$ of the graph of S.

CatTh[**FinLim**, S] is clearly closed under all the constructions of Appendix A. The properties listed in the preceding paragraph imply that the category CatTh[**FinLim**, S] is minimal with respect to the constructions of that Appendix, so that in fact those constructions can be taken as an recursive definition of CatTh[**FinLim**, S].

These observations are summed up in the following proposition.

**4.1.1 Proposition** *For a given sketch S, every object and every arrow of* CatTh[**FinLim**, S] *is constructible by repeated applications of the constructions of Appendix A.*

It is proved in [Barr and Wells, 1993a] that every object of CatTh[**FinLim**, S] is the limit of a diagram in the graph of the sketch S. Once this is known, the proof is a straightforward induction showing that the least fixed point of the constructions in Appendix A is indeed a category with finite limits.

## 4.2 Remark concerning models

Let *Cat* be some given finite-limit sketch for categories. We have to consider the following three mathematical entities.

1. A category 𝒞.

2. A model ℭ of the sketch *Cat* in the category of sets which "is" or "represents" 𝒞. ℭ is a morphism of finite-limit sketches from *Cat* to the underlying finite-limit sketch of **Set**. ℭ is determined uniquely by 𝒞 if one assumes that **Set** has canonical finite limits. It is determined up to natural isomorphism in any case.

3. The model CatTh[**FinLim**, ℭ] of CatTh[**FinLim**, *Cat*] induced by ℭ. This is a finite-limit-preserving functor from CatTh[**FinLim**, *Cat*] to **Set**.

For many category theorists, 𝒞, ℭ and CatTh[**FinLim**, ℭ] denote the "same thing", or perhaps 𝒞 and CatTh[**FinLim**, ℭ] are the "same thing" and ℭ consists of presentation data for $\mathscr{C}$. Other mathematicians would disagree, saying $\mathscr{C}$ is a category (presumably for them "category" has some meaning other than "model for the sketch for categories") and CatTh[**FinLim**, *Cat*] is a functor, so how could they be the same thing? In the sequel, we will distinguish these constructions typographically, but we do not take a stand on this difference of opinion. Note that in any case the following three categories are equivalent:

1. The category of small categories and functors.

2. **Mod**[*Cat*, **Set**], the category of models of the sketch *Cat* in the category **Set** with natural transformations as morphisms.

3. The category of finite-limit-preserving functors from CatTh[**FinLim**, *Cat*] to **Set** with natural transformations between them as morphisms.

## 5 Sketches in general

We provide here a definition of "sketch" based on [Wells, 1990] (some of the terminology has been changed.) This constitutes a generalization of the concept of sketch invented by Charles Ehresmann and described in [Bastiani and Ehresmann, 1972] or in [Barr and Wells, 1990]. This definition presupposes the concept of finite-limit sketch (see Sections 1.1 and 4).

### 5.1 Constructor spaces

#### 5.1.1 Assumptions

We will assume given a fixed finite-limit sketch *Cat* whose category of models is the category of small categories and functors. A specific such sketch is given explicitly in Appendix B.2. The presentation that follows has *Cat* as an implicit parameter.

**5.1.2 Notation** For any finite-limit sketch *S* and morphism $\eta : \textit{Cat} \to \textit{S}$ of sketches, the induced functor between categorial theories will be denoted by

$$\mathsf{CatTh}[\mathbf{FinLim}, \eta] : \mathsf{CatTh}[\mathbf{FinLim}, \textit{Cat}] \to \mathsf{CatTh}[\mathbf{FinLim}, \textit{S}]$$

**5.1.3 Definition** A finite-limit sketch $E$ together with a morphism $\eta : \mathsf{Cat} \to E$ of sketches is called a **constructor space sketch** provided that every object in $\mathsf{CatTh}[\mathbf{FinLim}, E]$ is the limit of a finite diagram whose nodes are of the form $\mathsf{CatTh}[\mathbf{FinLim}, \eta](\mathsf{n})$, where $\mathsf{n}$ is a node of $\mathsf{Cat}$.

**5.1.4 Remark** In all the examples in this article, $\mathsf{CatTh}[\mathbf{FinLim}, \eta]$ is injective on objects. We will often assume that it is an inclusion and write $\mathsf{n}$ for $\mathsf{CatTh}[\mathbf{FinLim}, \eta](\mathsf{n})$.

**5.1.5 Definition** A category of the form $\mathsf{CatTh}[\mathbf{FinLim}, E]$ for some constructor space sketch $E$ is called a **constructor space**. We will normally denote the constructor space $\mathsf{CatTh}[\mathbf{FinLim}, E]$ by $\mathbf{E}$.

**5.1.6 Definition** A model of a constructor space $\mathbf{E}$ in **Set** is called an **E-category**, and a morphism of such models is called an **E-functor**.

**5.1.7 Discussion**

Definitions 5.1.3, 5.1.5 and 5.1.6 generalize those in [Wells, 1990]. They conform to the two-dimensional version given in [Power and Wells, 1992].

The value in a model $\mathfrak{M}$ of an object $\mathsf{v}$ in a constructor space is the set of all examples of a particular construction that is possible in the **E**-category $\mathfrak{M}$. Extended examples of this are discussed in Section 5.2 and Section 7. Thus each object of $\mathbf{E}$ represents a *type of construction* possible in $\mathbf{E}$, hence the name "constructor space".

Examples of constructor space sketches include the sketches for categories with finite products, categories with finite limits, cartesian closed categories, toposes, and so on. Several constructor space sketches are given in detail in Appendix B.

## 5.2 Notation for diagrams in a constructor space

The object of $\mathsf{CatTh}[\mathbf{FinLim}, \mathsf{Cat}]$ whose value in a model is the set of all not necessarily commutative diagrams of the form

$$\begin{array}{ccc} A & \xrightarrow{h} & B \\ f \downarrow & \overset{x}{\nearrow} & \downarrow k \\ C & \xrightarrow{g} & D \end{array} \tag{13}$$

is the limit of the diagram

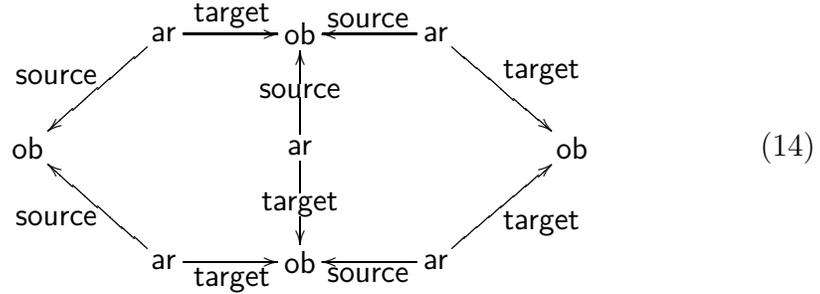
(14)

We now describe this in more detail and introduce some notation that makes the discussion of such diagrams easier to follow. We use the notation $D(n)$ to refer to the diagram shown herein with label $(n)$, and $I(n)$ for its shape graph. For example, the limit of the diagram above is $\mathsf{Lim}[D(14)]$.

Every node of $D(14)$ is either the object ob (the object that becomes the set of objects in a model) or the object ar (the object that becomes the set of arrows in a model) of CatTh[**FinLim**, *Cat*]. For a model $\mathfrak{C}$ of CatTh[**FinLim**, *Cat*] in **Set**, an element of $\mathfrak{C}(\mathsf{Lim}[D(14)]$ is a diagram in $\mathfrak{C}$, not necessarily commutative, of the form of Diagram (13).

In order to make the relation between Diagrams (13) and (14) clear, we give the shape graph of (14):

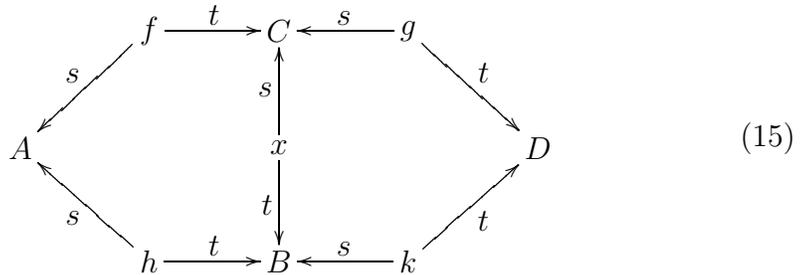
(15)

We have named the nodes of Diagram (15) by the objects and arrows that occur in Diagram (13) in such a way that the node named by an object or arrow of Diagram (13) will inhabit the value of that node in the model $\mathfrak{C}$. For example, the object $A$ of $\mathfrak{C}$ is at the upper left corner of Diagram (13) and the projection arrow from $\mathfrak{C}(\mathsf{Lim}[D(14)])$ to $\mathfrak{C}(\mathsf{ob})$ determined by the node labeled $A$ of the shape graph (15) is a function from the set of diagrams in $\mathfrak{C}$ of the form of Diagram (13) to the set of objects of $\mathfrak{C}$ that takes a diagram to the object in its upper left corner.

We will combine diagrams such as Diagram (14) and their shape graphs into one graph by labeling the nodes of the diagram by superscripts naming the corresponding node of the shape graph. In the case of Diagram (14), doing this gives the following **annotated diagram**:

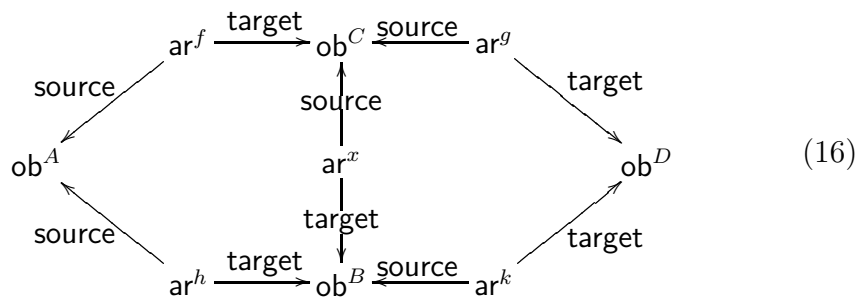 (16)

Here, the superscript $A$ on the leftmost node indicates that the corresponding node of the shape graph is labeled $A$. Formally, the expression $\mathsf{ob}^A$ is used as the label for the node $\delta(A)$ and its use signifies that $\delta(A) = \mathsf{ob}$. That device helps the reader to see that Diagram (13) is indeed an element of $\mathfrak{C}(\mathsf{Lim}[D(14)])$.

For example, the particular arrow $h$ of Diagram (13) is an element of $\mathsf{ar}$, and the label $\mathsf{ar}^h$ in Diagram (16) helps one see that it is that node that projects to $h$ in the model $\mathfrak{C}$ and that the source of $h$ is $A$ and that the target is $B$.

*It is important to understand that an annotated diagram such as (16) denotes precisely the same diagram as (14).* The fact that one node is labeled $\mathsf{ob}^A$ and another $\mathsf{ob}^B$ does not change the fact that *both nodes* are $\mathsf{ob}$. The superscript merely gives information about the relation between Diagram (14) and Diagram (13).

Diagram (16) could also be drawn as the base of a limit cone $\Theta$ with limit

Lim$[D(16)]$ (which of course is the same as Lim$[D(14)]$) as follows.

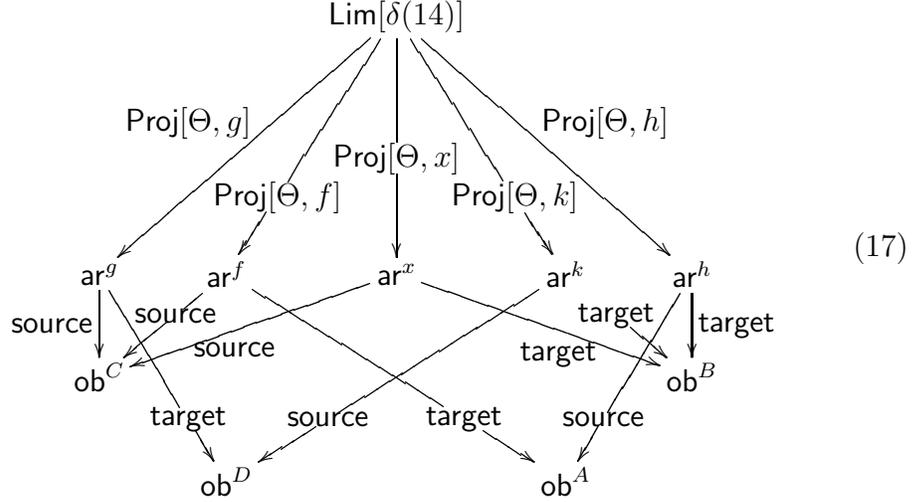

(17)

Because of the typographical complexity of doing this for diagrams more complicated than Diagram (16), we will usually give diagrams whose limits we discuss in the form of Diagram (16), without showing the cone, instead of in the form of Diagram (17).

Showing the cone explicitly as in Diagram (17) nevertheless has an advantage. It makes it clear that many of the projection arrows from Lim$[D(16)]$ are induced by others; in the particular case of Diagram (17), all the arrows to nodes of the form ob are induced by composing arrows to some node of the form ar with source or target. Diagram (16) does not exhibit this property. A systematic method of translating from graphical expressions such as Diagram (17) to a string-based expression could presumably be based on this. In the case of Diagram (17), the string-based expression would be something like this:

$$\{\langle g, f, x, k, h\rangle \mid \text{source}(g) = C, \text{target}(g) = D, \quad \text{source}(f) = A,$$
$$\text{target}(f) = C, \text{source}(x) = C, \text{target}(x) = B, \text{source}(k) = B,$$
$$\text{target}(k) = D, \text{source}(h) = A, \text{target}(h) = B\}$$

or in more familiar terms,

$$\{\langle g, f, x, k, h\rangle \mid g: C \to D, f: A \to C, x: C \to B, k: B \to D, h: A \to B\}$$

### 5.3 Sketches

**5.3.1 Definition** Let **E** be a constructor space. For any node v of **E**, we may freely adjoin a global element $S:1 \to v$ to obtain a finite-limit theory **E**[$S$]. $S$ is called an **E-sketch** and **E**[$S$] is called the **syntactic category** of $S$; we will denote it SynCat[$S$]. It should be distinguished from the categorial theory CatTh[**E**, $S$], which is the initial model of **E**[$S$]. Details are in [Wells, 1990] or [Power and Wells, 1992].

**5.3.2 Remark** If **E** = **FinLim**, then a finite-limit sketch $S$ in the traditional sense (a graph with diagrams and cones) corresponds to a **FinLim**-sketch in the sense of Definition 5.3.1 that we denote by Form[$S$] (but see Remark 5.3.4). The traditional sketch $S$ is an element of the value in the initial model of a certain (not uniquely determined) node v of **FinLim** which is the limit of a generally large and complicated diagram (not uniquely determined) that specifies the graph, diagrams and cones of $S$. Example 5.3.3 illustrates this correspondence. These remarks apply to other types of Ehresmann sketches by replacing **FinLim** by the suitable constructor space.

**5.3.3 Example** As an example, consider the **FinLim**-sketch $S$ with graph

$$A \xrightarrow[u]{v} B \xrightarrow{f} C \tag{18}$$

one diagram

$$\begin{array}{c} B \xrightarrow{u} A \\ \mathsf{Id}[B] \downarrow \swarrow v \\ B \end{array} \tag{19}$$

and one cone

$$\begin{array}{c} B \\ u \swarrow \searrow f \\ A \qquad C \end{array} \tag{20}$$

This is in fact a finite-product sketch, but any such sketch is also a **FinLim**-sketch. This is achieved, for instance, by setting $S := \mathsf{Form}[S]:1 \to v$, where v is the limit of this diagram:

$$\begin{array}{c}
\text{ob}^A \times \text{ob}^C \\
\downarrow \text{prod} \\
\text{ar}^u \xleftarrow{\text{lproj}} \text{cone} \xrightarrow{\text{rproj}} \text{ar}^f \xrightarrow{\text{target}} \text{ob}^C \\
\end{array} \tag{21}$$

(diagram with nodes $\text{ob}^A \times \text{ob}^C$, cone, $\text{ar}^u$, $\text{ar}^f$, $\text{ob}^A$, $\text{ob}^B$, $\text{ob}^C$, $\text{ar}_2$, ar, $\text{ar}^v$, with arrows: prod, lproj, rproj, target, source, target, rfac, source, comp, unit, lfac, target, source)

This is however not the only way of viewing $S$ as a **FinLim**-sketch as explained below.

**5.3.4 Remark** Note that the identification between traditional sketches and sketches in the new sense is many-to-one in both directions. A traditional finite-limit sketch can be described as an element of the limit of more than one diagram in **FinLim**. In the other direction, because every node v of **FinLim** is the limit of a diagram in a sketch for **FinLim** (for example, the one in Appendix B.5), any element of the value of v in an initial model can be described in terms of a graph, diagrams and cones, but in general not in only one way, and in any case such a description requires an arbitrary choice for the names of the nodes and arrows of the graph.

### 5.4 Models of sketches

Let $S$ be an **E**-sketch as in Definition 5.3.1. The category of models (finite-limit-preserving functors) of **E**$[S]$ has an underlying functor $\mathsf{UndCat}[\mathbf{E}, S]$ to the category of **E**-categories induced by the embedding of **E** into **E**$[S]$.

**5.4.1 Definition** Let $S$ be an **E**-sketch and $\mathscr{C}$ an **E**-category. An **E**-functor from $\mathsf{CatTh}[\mathbf{E}, S]$ to $\mathscr{C}$ is called a **model** of $S$. A **morphism of models** in $\mathscr{C}$ is a natural transformation. The category of models of $S$ in $\mathscr{C}$ is denoted by $\mathbf{Mod}[S, \mathscr{C}]$.

**5.4.2 Remark** Each model of an **E**-sketch $S$ is the unique functor induced (because $\mathsf{CatTh}[\mathbf{E}, S]$ is an initial object) by a model $\mathfrak{M}$ of **E**$[S]$ in **Set** (a finite-limit-preserving functor from $\mathsf{E}[S]$ to **Set**) that has $\mathscr{C}$ as its underlying **E**-category. This passage from model of $S$ to model of $\mathsf{E}[S]$ induces an equivalence of categories between the category of **E**-functors from $\mathsf{CatTh}[\mathsf{E}, S]$ to $\mathscr{C}$ and models (**FinLim**-functors) $\mathfrak{M}$ of $\mathsf{E}[S]$ in **Set** that have $\mathscr{C}$ as underlying **E**-category.

**5.4.3 Example** In Example 5.3.3, a model of the traditional sketch $S$ in a finite-limit category $\mathscr{C}$ in the sense of Definition 5.4.1 is a finite-limit preserving functor from $\mathsf{CatTh}[\mathsf{FinLim}, S]$ (which is the initial model of $\mathsf{FinLim}[S]$ that takes $\mathsf{S}$ to $S$) to $\mathscr{C}$. The finite-limit functor from $\mathsf{FinLim}[S]$ to **Set** mentioned in Remark 5.4.2 is the extension of the model $\mathfrak{C}\colon \mathsf{FinLim} \to \mathbf{Set}$ that "corresponds to" or "is" (see Remark 4.2) $\mathscr{C}$ obtained by setting the value of the global element $\mathsf{S}$ to be the traditional sketch $S$.

# 6 Assertions in graph-based logic

## 6.1 Potential factorizations

In a constructor space **E**, a diagram of the form

$$\begin{array}{c} C \\ \downarrow j \\ A \xrightarrow{f} B \end{array} \tag{22}$$

is called a **potential factorization** or **PF**. $B$ is called the **workspace** of the potential factorization, and $j$ is called the **claim** of the potential factorization. In many examples, including all theories in this paper, $j$ is monic and corresponds to a formal selection of a subset of those objects formally denoted by $B$.

### 6.1.1 Notation

We will use suggestive notation for a potential factorization, as exhibited in the diagram below.

$$\begin{array}{c} \textit{claim} \\ \downarrow \textit{claimcon} \\ \textit{hyp} \xrightarrow[\textit{hypcon}]{} \textit{wksp} \end{array} \tag{23}$$

The names *hypcon*, *claimcon* and *wksp* respectively abbreviate "hypothesis construction", "claim construction" and "workspace". The reason for the names of the arrows and nodes is discussed in 7.3.

### 6.1.2 Actual factorizations

If an actual factorization arrow *verif* : *hyp* → *claim* exists in **E** that makes

$$\begin{array}{c} & claim \\ \nearrow^{verif} & \downarrow^{claimcon} \\ hyp \xrightarrow[hypcon]{} wksp \end{array} \qquad (24)$$

commute, then we say that the potential factorization is **deducible**. If for some model $\mathfrak{M}$ of **E** there is an arrow $\xi$ of **Set** that makes

$$\begin{array}{c} & \mathfrak{M}(claim) \\ \nearrow^{\xi} & \downarrow^{\mathfrak{M}(claimcon)} \\ \mathfrak{M}(hyp) \xrightarrow[\mathfrak{M}(hypcon)]{} \mathfrak{M}(wksp) \end{array} \qquad (25)$$

commute, then we say the model $\mathfrak{M}$ **satisfies** the potential factorization. If for *every* model $\mathfrak{M}$ there is such an arrow $\xi$ then we say that the potential factorization is **valid**.

We give examples of potential factorizations in Section 7, and then discuss the general concept of potential factorization in Section 7.3.

## 6.2 Completeness

**6.2.1 Theorem** *In any categorial theory* CatTh[**E**, *S*], *a potential factorization is deducible if and only if it is valid.*

**Proof** That deducibility implies validity follows from the fact that functors preserve factorizations.

For the converse, let

$$\begin{array}{c} claim \\ \downarrow^{claimcon} \\ hyp \xrightarrow[hypcon]{} wksp \end{array} \qquad (26)$$

be a potential factorization in $\mathsf{CatTh}[\mathbf{E}, S]$. Suppose it is valid. Because the functor $\mathrm{Hom}(\textit{hyp}, -)$ is a model, the hypothesis of the theorem implies that we may choose an arrow $\xi$ of **Set** such that the diagram

$$\begin{array}{c}
\mathrm{Hom}(\textit{hyp}, \textit{claim}) \xrightarrow{\xi} \mathrm{Hom}(\textit{hyp}, \textit{claim}) \\
\downarrow \mathrm{Hom}(\textit{hyp}, \textit{claimcon}) \\
\mathrm{Hom}(\textit{hyp}, \textit{claim}) \xrightarrow{\mathrm{Hom}(\textit{hyp}, \textit{hypcon})} \mathrm{Hom}(\textit{hyp}, \textit{wksp})
\end{array} \quad (27)$$

commutes. Define $\textit{verif} := \xi(\mathsf{id}_{\textit{hyp}})$. Then $\textit{claimcon} \circ \textit{verif} = \textit{hypcon}$, so that the potential factorization is deducible.

Any model $\mathfrak{M}$ is a functor that preserves finite limits, so that if $\textit{claimcon}$ is monic, so is $\mathfrak{M}(\textit{claimcon})$. In that case, $\mathfrak{M}(\textit{verif})$ necessarily equals $\xi$. □

Section B.7.2 gives an example of how completeness can be used.

**6.2.2 Remark** Traditional proofs of completeness contain an induction which is missing in the preceding argument. In the present system, a theorem can be identified with an arrow of $\mathsf{CatTh}[\mathbf{FinLim}, S]$. Proposition 4.1.1 describes the recursive construction of arrows and so is the analog of the inductive part of traditional proofs of completeness.

**6.2.3 Proposition** *Let $D : I \to \mathsf{CatTh}[\mathbf{E}, S]$ be a diagram. Suppose for every model $\mathfrak{M}$ of $\mathsf{CatTh}[\mathbf{E}, S]$, $\mathfrak{M} \circ D$ commutes. Then $D$ commutes.*

**Proof** Suppose

$$\begin{array}{c}
A \xrightarrow{f} B \\
{}_{h} \searrow \quad \downarrow g \\
C
\end{array} \quad (28)$$

is a diagram in $\mathsf{CatTh}[\mathbf{E}, S]$. Then, because $\mathrm{Hom}(A, -)$ is a model,

$$\begin{array}{c}
\mathrm{Hom}(A, A) \xrightarrow{\mathrm{Hom}(A, f)} \mathrm{Hom}(A, B) \\
{}_{\mathrm{Hom}(A, h)} \searrow \quad \downarrow \mathrm{Hom}(A, g) \\
\mathrm{Hom}(A, C)
\end{array}$$

commutes. By chasing Id[A] around the diagram both ways, we get $g \circ f = h$, so Diagram (28) commutes as well.□

## 7 Examples of potential factorizations

We give two examples of potential factorizations in some detail. The example in 7.1 holds for any category and the example in 7.2 holds in any Cartesian closed category. We construct actual factorizations for these PF's in Section 9.

### 7.1 A fact about diagrams in any category

The following proposition holds in any category.

**7.1.1 Proposition** *In any category, given the following diagram*

$$
\begin{array}{ccc}
A & \xrightarrow{h} & B \\
f \downarrow & \overset{x}{\underset{g}{\diagup}} & \downarrow k \\
C & \xrightarrow{g} & D
\end{array}
\tag{29}
$$

*if the two triangles commute then so does the outside square.*

**Proof**  $k \circ h = k \circ (x \circ f) = (k \circ x) \circ f = g \circ f.$ □

We construct here the potential factorization in CatTh[**FinLim**, *Cat*] that corresponds to this theorem.

As we pointed out in Section 5.2, the value in a model of Lim[D(16)] is a diagram such as Diagram (29). However, Diagram (16) only carries the information as to the source and target of the arrows in Diagram (29). The structure we must actually work with includes the information as to which pairs are composable. There are four composable pairs in Diagram (29), and each one inhabits the value of $ar_2$, the node of formal composable pairs in a category (see Appendix B.2).

The following diagram thus contains the basic information about sources, targets and composability that are required for stating the theorem, and so its limit is suitable to be *wksp*.

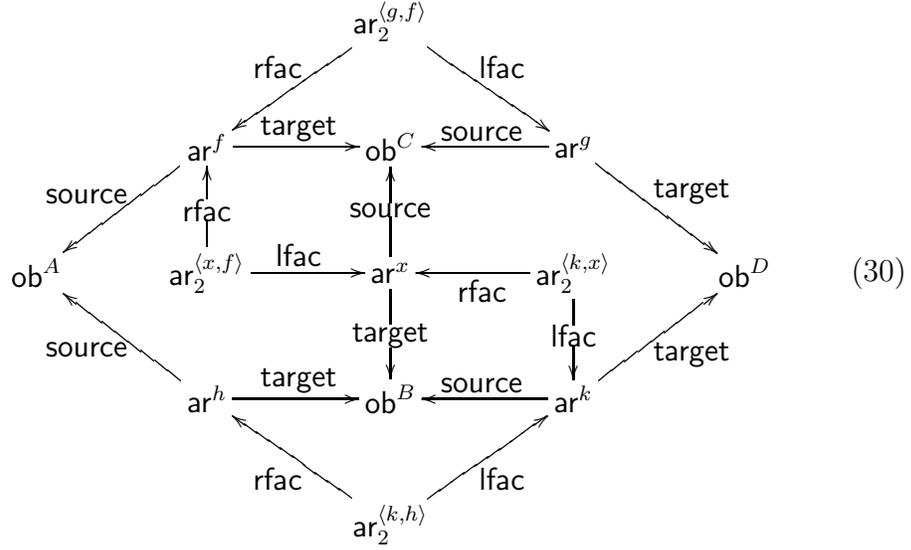

(30)

The statement that the two triangles in Diagram (29) commute is: $x \circ f = h$ and $k \circ x = g$. Using the composition arrow $\mathsf{comp}:\mathsf{ar}_2 \to \mathsf{ar}$ of $\mathsf{CatTh}[\mathbf{FinLim}, \mathit{Cat}]$, this statement amounts to saying that Diagram (29) is a member of $\mathfrak{C}(\mathsf{Lim}[D(31)])$, which we take as *hyp*:

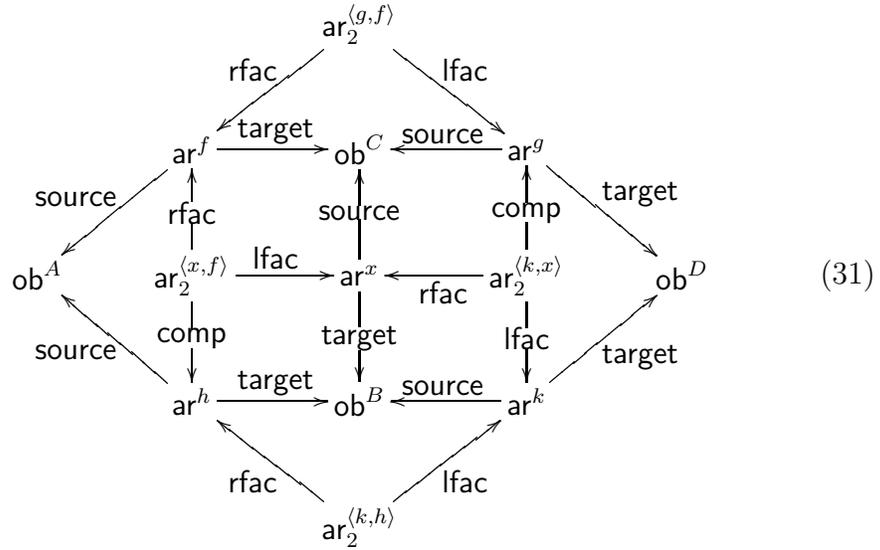

(31)

The statement that the outside of Diagram (29) commutes is that the diagram is a member of $\mathfrak{C}(\mathsf{Lim}[\delta(32)])$, which we take as *claim*:

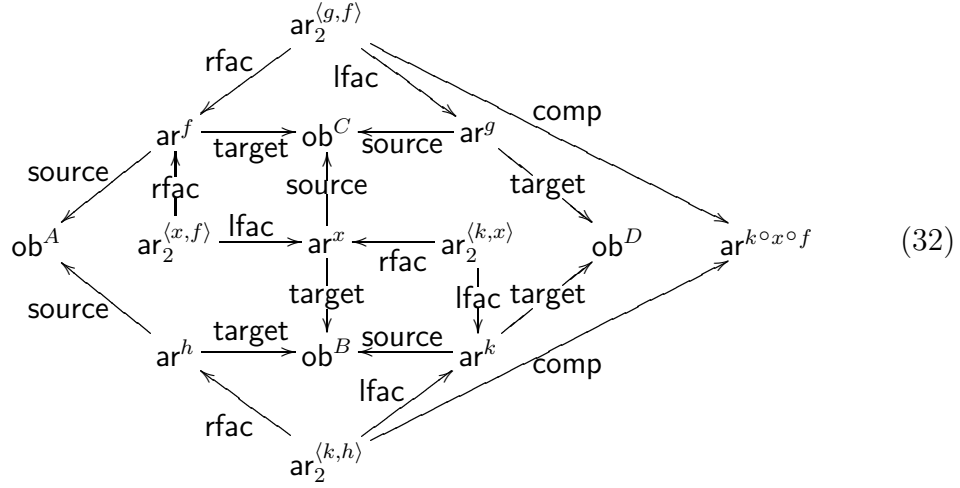
(32)

Because Diagram (30) and Diagram (31) are restrictions of Diagram (32), we may, by Lemma 3.2.1 choose a function

$$\mathit{claimcon} : \mathsf{Lim}[D(32)] \to \mathsf{Lim}[D(30)]$$

and similarly another function

$$\mathit{hypcon} : \mathsf{Lim}[D(31)] \to \mathsf{Lim}[D(30)]$$

producing a potential factorization

$$\begin{array}{c} & \mathsf{Lim}[D(32)] \\ & \downarrow \mathit{claimcon} \\ \mathsf{Lim}[D(31)] \xrightarrow[\mathit{hypcon}]{} \mathsf{Lim}[D(30)] \end{array} \quad (33)$$

This potential factorization expresses the content of Proposition 7.1.1 in diagrammatic form. It should be clear that the node $\mathsf{Lim}[D(30)]$ could have been replaced by $\mathsf{Lim}[D(16)]$.

## 7.2 A fact about Cartesian closed categories

Proposition 7.1.1 holds in any category. We now discuss a theorem of Cartesian closed categories, to show how the system presented in this paper handles structure that cannot be expressed using Ehresmann sketches. The latter are equivalent in expressive power to ordinary first order logic. (An excellent presentation of the details of this fact may be found in [Adámek and Rosičky, 1994].) Thus this example in a certain sense requires higher-order logic.

### 7.2.1 Proposition  *In any Cartesian closed category, if*

$$
\begin{array}{ccc}
A \times B & \xrightarrow{g} & C \\
& \searrow{h} & \downarrow{f} \\
& & D
\end{array}
$$

*commutes, then so does*

$$
\begin{array}{ccc}
A & \xrightarrow{\lambda g} & C^B \\
& \searrow{\lambda h} & \downarrow{f^B} \\
& & D^B
\end{array}
$$

See Appendix B.6 for notation. The arrow $f^B$ is defined by

$$f^B := \lambda(C^B \times B \xrightarrow{\text{eval}} C \xrightarrow{f} D) : C^B \longrightarrow D^B \tag{34}$$

The proof follows from the fact that $\lambda$ is invertible and the calculation

$$
\begin{aligned}
\text{eval} \circ \left((f^B \circ \lambda g) \times \text{Id}[B]\right) &= \text{eval} \circ \left((f^B \times \text{Id}[B]) \circ (\lambda g \times \text{Id}[B])\right) \\
&= \left(\text{eval} \circ (f^B \times \text{Id}[B])\right) \circ (\lambda g \times \text{Id}[B]) \\
&= (f \circ \text{eval}) \circ (\lambda g \times \text{Id}[B]) \\
&= f \circ (\text{eval} \circ (\lambda g \times \text{Id}[B])) \\
&= f \circ g = h = \text{eval} \circ (\lambda h \times \text{Id}[B])
\end{aligned}
$$

The first equality is based on an assertion true in all categories that can be handled in our system in a manner similar to (but more complicated than) that of module (63) in Appendix B.4. The second and fourth equalities are associativity of composition and are proven using Figure (53) of Appendix B.2. The sixth equality is a hypothesis. The other three equalities are all based on Diagram (83) in Appendix B.6.

We present here the potential factorization corresponding to the third equality, which is the most complicated of those based on Diagram (83). In this presentation, unlike that of 7.1, we will use the modules developed in Appendix B to simplify the figures. The actual factorization corresponding to this potential factorization is given in 9.2.

The fact under discussion is that the diagram

$$
\begin{array}{ccccc}
A \times B & \xrightarrow{\lambda g \times \mathsf{Id}[B]} & C^B \times B & \xrightarrow{f^B \times \mathsf{Id}[B]} & D^B \times B \\
& & \downarrow{\scriptstyle \mathrm{eval}} & & \downarrow{\scriptstyle \mathrm{eval}} \\
& & C & \xrightarrow{f} & D
\end{array}
\tag{35}
$$

commutes.

Thus *wksp* will be the limit of the following diagram, which describes the objects and arrows in Diagram (35) but has no requirements on its commutativity.

$$
\begin{array}{ccccccc}
\mathsf{ar}^{\lambda g \times \mathsf{Id}[B]} & \xrightarrow{\mathsf{target}} & \mathsf{ob}^{C^B \times \mathsf{Id}[B]} & \xleftarrow{\mathsf{source}} & \mathsf{ar}^{f^B \times \mathsf{Id}[B]} & \xrightarrow{\mathsf{target}} & \mathsf{ob}^{D^B \times \mathsf{Id}[B]} \\
\downarrow{\scriptstyle \mathsf{source}} & & \uparrow{\scriptstyle \mathsf{source}} & & & & \uparrow{\scriptstyle \mathsf{source}} \\
\mathsf{ob}^{A \times B} & & \mathsf{ar}^{\mathrm{eval}} & & & & \mathsf{ar}^{\mathrm{eval}} \\
& & \downarrow{\scriptstyle \mathsf{target}} & & & & \downarrow{\scriptstyle \mathsf{target}} \\
& & \mathsf{ob}^C & \xleftarrow{\mathsf{source}} & \mathsf{ar}^f & \xrightarrow{\mathsf{target}} & \mathsf{ob}^D
\end{array}
\tag{36}
$$

We define *hyp* to be the limit of the following diagram, in which $\phi =$

$(f \circ \mathsf{eval}) \circ (\lambda g \times \mathsf{Id}[B])$.

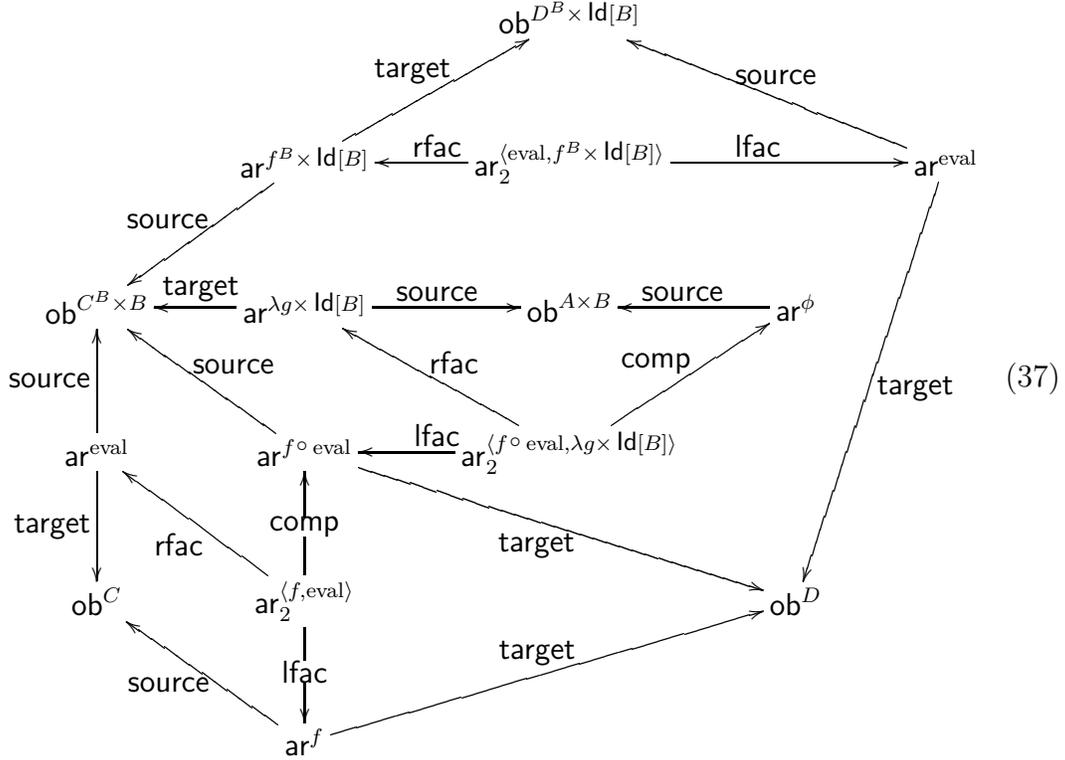

(37)

Then Diagram (36) is a subdiagram of Diagram (37) (the big rectangle in Diagram (36) is the perimeter of Diagram (37)) and we define *hypcon* to be the induced arrow from $\mathsf{Lim}[D(37)]$ to $\mathsf{Lim}[D(36)]$.

*claim* is the limit of a diagram we shall refer to as Diagram (37′), obtained from Diagram (37) by adjoining an arrow labeled comp from $\mathsf{ar}_2^{\langle \mathsf{eval}, f^B \times B \rangle}$ to $\mathsf{ar}^{f \circ \mathsf{eval}}$. Diagram (37′) includes Diagram (37) and hence Diagram (36), and we take the arrow *claimcon* to be the arrow from Diagram (37′) to Diagram (36) induced by this inclusion. The existence of both *claimcon* and *hypcon* follows from Lemma 3.2.1.

We now have a potential factorization

$$\begin{array}{c} \mathsf{Lim}[D(37')] \\ \downarrow {\scriptstyle claimcon} \\ \mathsf{Lim}[D(37)] \xrightarrow[hypcon]{} \mathsf{Lim}[D(36)] \end{array} \qquad (38)$$

The above potential factorization expresses the content of the third equality in the calculation in the proof of Proposition 7.2.1.

## 7.3 Discussion of the examples

### 7.3.1 General discussion

In a potential factorization

$$\begin{array}{c} \textit{claim} \\ \downarrow \textit{claimcon} \\ \textit{hyp} \xrightarrow[\textit{hypcon}]{} \textit{wksp} \end{array} \qquad (39)$$

each of *hyp*, *claim* and *wksp* represents a type of entity that can be constructed in an **E**-category. In this section, we discuss a way of thinking about these nodes that exhibits how they could represent a theorem about **E**-categories. This discussion is relevant to potential factorizations that are motivated by theorems stated in the traditional manner; it may not be relevant to an arbitrary potential factorization, which after all can be any pair of arrows with a common codomain.

1. The node *wksp* (for "workspace") represents the data involved in both the hypothesis and the conclusion. For a given theorem, the choice of what is actually included in *wksp* may be somewhat arbitrary (see the discussion of Example 7.2 below).

2. The node *hyp* represents possible additional properties that are part of the assumptions in the theorem being represented.

3. The node *claim* represents the properties that the theorem asserts must hold given the assumptions.

4. The arrow *hypcon* represents the selection or construction necessary to see the hypothesis as part of the workspace. In both our examples, *hypcon* represents a simple forgetting of properties.

5. The arrow *claimcon* represents the selection or construction necessary to see the claim as part of the workspace.

6. The arrow *verif* in an actual factorization represents a specific way, *uniform in any model*, that any entity of type *hyp* can be transformed into, or recognized as, an entity of type *claim*.

### 7.3.2 Discussion of Example 7.1

In Example 7.1,

1. *wksp* represent squares of the form

$$A \xrightarrow{h} B \quad f \downarrow \ \overset{x}{\nearrow} \ \downarrow k \quad C \xrightarrow{g} D \tag{40}$$

   with no commutativity conditions.

2. *hyp* represents diagrams of the form Diagram (40) in which the two triangles commute;

3. *claim* represents diagram of the form Diagram (40) in which the outside square commutes.

4. *hypcon* represents forgetting that the two triangles commute. Because it represents forgetting a property in this case, *hypcon* is monic, but in general it need not be.

5. *claimcon* represents forgetting that the outside square commutes.

### 7.3.3 Discussion of Example 7.2

In Example 7.2,

1. *wksp* represent diagrams of the form Diagram (35) with no commutativity conditions and no recognition that any sequence of arrows is composable. The phrase "The form of Diagram (35)" refers to the source and target commonalities of the arrows in the diagram. Obviously, some sequences compose but we have not represented that in *wksp*, although we could have.

2. *hyp* represents diagrams of the form of Diagram (35), recognizing the composite $f \circ \mathrm{eval}$ and the fact that

$$f \circ \mathrm{eval}, \lambda g \times \mathsf{Id}[B]$$

and $\mathrm{eval}, f^B \times \mathsf{Id}[B]$ are composable pairs.

3. *claim* represents diagram of the form Diagram (35) that commute.

4. *hypcon* represents forgetting the information concerning composition in *hyp*.

5. *claimcon* represents forgetting the information concerning composition in *claim*.

We discuss the meaning of the actual factorization arrows *verif* in 9.3.

**7.3.4 Remark** The representation of facts such as those of Example 7.1 and 7.2 as potential factorizations is *variable free* in the sense that in each statement, one does not refer to a particular diagram such as Diagram (29) or Diagram (35) which stands as a pattern for all such diagrams. Propositions 7.1.1 and 7.2.1 state the fact in question using those diagrams as patterns, and understanding their meaning calls on the reader's ability to recognize patterns. Our description of the fact in the examples as a potential factorization is much more complicated because the diagrams involved in the potential factorization are essentially explicit descriptions of the relations between the nodes and arrows of Diagram (29) and Diagram (35) respectively, relations which a knowledgeable reader grasps from seeing the diagrams without having them indicated explicitly.

# 8 Construction and Deduction

## 8.1 The rules of graph-based logic

In first-order logic, rules are given for constructing terms and formulas, and further rules (rules of inference) are given for deriving formulas from formulas. These rules are intended to preserve truth.

In Appendix A we give rules for constructing all the objects and arrows of CatTh[**FinLim**, $S$] for an arbitrary finite-limit sketch $S$, and for constructing a basis for all the commutative diagrams in CatTh[**FinLim**, $S$]. These rules correspond to both the term and formula construction rules and the rules of inference of string-based logic. The tools of a typical string-based logic include constant symbols, variables, function symbols, logical operators and quantifiers. Here we have nodes, arrows and commutative diagrams. What corresponds to a sentence is a pair of arrows as in Diagram (39) (a potential factorization), and what corresponds to the satisfiability of the sentence in a model $\mathfrak{M}$ is the existence of an arrow $\xi$ for which

$$\begin{array}{ccc}
& & \mathfrak{M}(hyp) \\
& \xi \nearrow & \downarrow \mathfrak{M}(claimcon) \\
\mathfrak{M}(claim) & \xrightarrow[\mathfrak{M}(hypcon)]{} & \mathfrak{M}(wksp)
\end{array} \qquad (41)$$

commutes. A demonstration of the deducibility of the sentence corresponds to the construction of an arrow $verif : claim \to hyp$ in $\mathsf{CatTh}[\mathbf{FinLim}, S]$ such that

$$\begin{array}{ccc}
& & hyp \\
& verif \nearrow & \downarrow claimcon \\
claim & \xrightarrow[hypcon]{} & wksp
\end{array} \qquad (42)$$

commutes. Thus the same rules suffice for constructing the sentence (the potential factorization) and for proving it (constructing the arrow that makes it an actual factorizaton). (See Section 7.3 for further discussion of these points.)

Each rule in Appendix A is actually a rule scheme, and each instance of the scheme is an assertion that, given certain arrows and commutative diagrams (see Remark 8.1.1 below) in $\mathsf{CatTh}[\mathbf{FinLim}, S]$, other arrows or commutative diagrams exist in $\mathsf{CatTh}[\mathbf{FinLim}, S]$. For example, if $\delta$ is a diagram with shape graph

$$\begin{array}{c}
i \\
\downarrow u \\
j \xrightarrow{v} k
\end{array} \qquad (43)$$

then the following rule is an instance of ∃LIM:

$$\frac{\begin{array}{c}\delta(j)\xrightarrow{\delta(u)}\delta(k)\\ \phantom{x}\\ \delta(i)\downarrow\delta(u)\\ \phantom{x}\end{array}}{\mathsf{Lim}[\delta]\xrightarrow{\mathsf{Proj}[\mathsf{Lim}[\delta],i]}\delta(i),\ \mathsf{Proj}[\mathsf{Lim}[\delta],j],\ \mathsf{Proj}[\mathsf{Lim}[\delta],k]} \tag{44}$$

### 8.1.1 Remark

Rules ∃FIA, !FIA and CFIA assume the existence of commutative cones, but a commutative cone is a collection of interrelated commutative diagrams, so the statement above that each scheme assumes the existence of certain arrows and commutative diagrams is correct. Thus given a cone $\Theta$, an instance of ∃FIA is this rule:

$$\frac{\mathsf{Vertex}[\Theta]\xrightarrow{\mathsf{Proj}[\Theta,i]}\delta(i),\ \mathsf{Proj}[\Theta,j],\ \mathsf{Proj}[\Theta,k]}{\mathsf{Vertex}[\Theta]\xrightarrow{\mathsf{Fillin}[\Theta,\delta]}\mathsf{Lim}[\delta]} \tag{45}$$

The point of this remark is that ∃FIA is a rule with a diagram as hypothesis and an arrow as conclusion. The hypothesis is the cone itself, not the string "$\Theta : v \lhd (\delta : I \to \mathscr{C})$" or any other description of it.

## 9 Examples of theorems

### 9.1 Proof of Proposition 7.1.1

We continue Example 7.1 by constructing and thereby deducing the existence of an arrow *verif* : $\mathsf{Lim}[D(31)] \to \mathsf{Lim}[D(32)]$ making Diagram (33) commute.

We first construct Diagram (31′) (not shown) by adjoining $\mathsf{ar}_3$ to Diagram (31), along with arrows

$\mathsf{lfac} : \mathsf{ar}_3 \to \mathsf{ar}^k$
$\mathsf{mfac} : \mathsf{ar}_3 \to \mathsf{ar}^x$
$\mathsf{rfac} : \mathsf{ar}_3 \to \mathsf{ar}^f$

We further construct Diagram (31″) by adjoining

$\langle \mathsf{lfac}, \mathsf{mfac} \rangle : \mathsf{ar}_3 \to \mathsf{ar}_2^{\langle k,x \rangle}$
$\langle \mathsf{mfac}, \mathsf{rfac} \rangle : \mathsf{ar}_3 \to \mathsf{ar}_2^{\langle x,f \rangle}$
$\mathsf{lass} : \mathsf{ar}_3 \to \mathsf{ar}^{\langle g,f \rangle}$
$\mathsf{rass} : \mathsf{ar}_3 \to \mathsf{ar}^{\langle k,h \rangle}$

The first two are induced by the fact that $\mathsf{ar}_3$ is a limit and the other two are defined in Appendix B.2.

Diagram (30) is a base restriction of each of Diagram (31′) and Diagram (31″), so, we may, using Lemma 3.2.1, choose arrows

$$\phi_1 : \mathsf{Lim}[D(31')] \to \mathsf{Lim}[D(30)]$$

and

$$\phi_2 : \mathsf{Lim}[D(31'')] \to \mathsf{Lim}[D(30)]$$

Diagram (31) is a dominant subdiagram of Diagram (31′) since the latter is obtained from the former by adjoining a limit of a subdiagram together with their projection arrows ($\mathsf{lfac}$, $\mathsf{mfac}$ and $\mathsf{rfac}$). Therefore, using Lemma 3.3.4, we may choose an isomorphism $\psi_1$ making

$$\begin{array}{ccc} \mathsf{Lim}[D(31)] & \xrightarrow{\psi_1} & \mathsf{Lim}[D(31')] \\ & \searrow^{\mathsf{hypcon}} \quad \swarrow_{\phi_1} & \\ & \mathsf{Lim}[D(30)] & \end{array} \qquad (46)$$

commute.

Similarly, Diagram (31′) is a dominant subdiagram of Diagram (31″) since the latter is obtained from the former by adjoining four composites. Therefore by Lemma 3.3.4 we may choose an isomorphism $\psi_2$ making

$$\begin{array}{ccc} \mathsf{Lim}[D(31')] & \xrightarrow{\psi_2} & \mathsf{Lim}[D(31'')] \\ & \searrow^{\phi_1} \quad \swarrow_{\phi_2} & \\ & \mathsf{Lim}[D(30)] & \end{array} \qquad (47)$$

commute. We then construct Diagram $(31''')$ by adjoining arrows

$$\mathsf{comp} : \mathsf{ar}_2^{\langle g,f \rangle} \to \mathsf{ar}^{k \circ x \circ f}$$

and

$$\mathsf{comp} : \mathsf{ar}_2^{\langle k,h \rangle} \to \mathsf{ar}^{k \circ x \circ f}$$

where $\mathsf{ar}^{k \circ x \circ f}$ is a new node.

Because of associativity (the right diagram in Figure (53) of Appendix B.2), Diagram $(31''')$ extends Diagram $(31'')$ by adjoining composites, so we may choose an isomorphism $\psi_3 : \mathsf{Lim}[D(31'')] \to \mathsf{Lim}[D(31''')]$ and an arrow $\phi_3 : \mathsf{Lim}[D(31''')] \to \mathsf{Lim}[D(30)]$ making

$$\begin{array}{c}
\mathsf{Lim}[D(31'')] \xrightarrow{\psi_3} \mathsf{Lim}[D(31''')] \\
\phi_2 \searrow \quad \swarrow \phi_3 \\
\mathsf{Lim}[D(30)]
\end{array} \qquad (48)$$

commute.

Finally, by Lemma 3.2.1, we may choose arrow $\psi_4 : \mathsf{Lim}[D(31''')] \to \mathsf{Lim}[D(32)]$ making

$$\begin{array}{c}
\mathsf{Lim}[D(31''')] \xrightarrow{\psi_4} \mathsf{Lim}[D(32)] \\
\phi_3 \searrow \quad \swarrow \mathsf{claim} \\
\mathsf{Lim}[D(30)]
\end{array} \qquad (49)$$

commute.

We next set

$$\mathit{verif} := \psi_4 \circ \psi_3 \circ \psi_2 \circ \psi_1 \qquad (50)$$

whence the theorem follows.

### 9.2 Proof of Theorem 7.2.1

In this section, we provide a factorization of the potential factorization described in Section 7.2.

Diagram (37) contains the following as a subdiagram, in which, using Diagram (34),

$$\theta = \langle \mathsf{eval}, f^B \times \mathsf{Id}[B] \rangle = \langle \mathsf{eval}, \lambda(f \circ \mathsf{eval}) \times \mathsf{Id}[B] \rangle$$

$$\begin{array}{c}
\text{ob}^{C^B \times B} \xleftarrow{\text{source}} \text{ar}^{f^B \times \text{Id}[B]} \xrightarrow{\text{target}} \text{ob}^{D^B \times B} \\
\uparrow \text{rfac} \\
\text{source} \uparrow \quad \text{ar}_2^{\theta} \quad \uparrow \text{source} \\
\swarrow \text{comp} \quad \text{lfac} \searrow \\
\text{ar}^{f \circ \text{eval}} \xrightarrow{\text{target}} \text{ob}^D \xleftarrow{\text{target}} \text{ar}^{\text{eval}}
\end{array} \qquad (51)$$

This diagram is an instance of Diagram (83), so it commutes. The arrow comp satisfies Definition 3.3.1, so Lemma 3.3.4 implies that there is an isomorphism *verif* : *hyp* → *claim*. Now the inclusion of Diagram (36) into Diagram (37) followed by the inclusion of Diagram (37) into Diagram (37′) is precisely the inclusion of Diagram (36) into Diagram (37′). It follows that *hypcon* ∘ *verif*$^{-1}$ = *claimcon*, so that *claimcon* ∘ *verif* = *hypcon* as required.

### 9.3 Discussion of the proofs.

The factorization *verif* : $\text{Lim}[D(31)] \to \text{Lim}[D(32)]$ of Diagram (33) given in Equation (50) constitutes the recognition that if the two triangles commute, then so does the outside square. The fact that *verif* makes Diagram (33) commute is a codification of the fact that if the two triangles commute then so does the outside square *of the same diagram.* In general, the reason we require that actual factorizations be an arrow in a certain comma category instead of merely an arrow from one node to another is to allow us to assert hypotheses and conclusions that share data (in this case the data in Diagram (40)).

The factorization *verif* : *hyp* → *claim* constructed in 9.2 constitutes recognition that $\text{eval} \circ (f^B \times \text{Id}[B]) = f \circ \text{eval}$ via the arrow $\text{comp} \colon \text{ar}_2^{\langle \text{eval}, f^B \times \text{Id}[B] \rangle} \to \text{ar}^{f \circ \text{eval}}$ in Diagram (37′). Because of this, the node $\text{ar}_2^{\langle f \circ \text{eval}, \lambda g \times \text{Id}[B] \rangle}$ could also be labeled $\langle \text{eval} \circ f^B \times \text{Id}[B], \lambda g \times \text{Id}[B] \rangle$. Thus the factorization also exhibits the fact that

$$\text{eval} \circ (f^B \times \text{Id}[B]) \circ (\lambda g \times \text{Id}[B]) = f \circ \text{eval} \circ (\lambda g \times \text{Id}[B])$$

It is clear that there are many alternative formulations of Proposition 7.2.1. For example, instead of first constructing $\text{ar}_2^{f \circ \text{eval}}$ as in Diagram (37) (which is *hyp* in this case), we could have constructed a node

$\mathsf{ar}_2^{\mathrm{eval}\,\circ\,(f^B\times\mathsf{Id}[B])}$ and an arrow

$$\mathsf{comp} : \mathsf{ar}_2^{\langle\mathrm{eval},(f^B\times\mathsf{Id}[B])\rangle} \to \mathsf{ar}_2^{\mathrm{eval}\,\circ\,(f^B\times\mathsf{Id}[B])}$$

Then the construction of an arrow

$$\mathsf{comp} : \mathsf{ar}_2^{\langle f,\mathrm{eval}\rangle} \to \mathsf{ar}^{\mathrm{eval}\,\circ\,(f^B\times\mathsf{Id}[B])}$$

would have proved the theorem.

## 10 Future work

A second part of this work [Bagchi and Wells, 1993a] will discuss in detail the relationship between the graph-based logic for finite-product categories and the rules for equational deduction.

It is noteworthy that the rules of construction for constructor spaces given in Appendix A correspond to arrows of **FinLim**, although not in a one-to-one way (see Remark A.1.1). The rules are given here in a form that requires pattern recognition (recall the discussion in 7.3.4), but they clearly could be given at another level of abstraction as arrows or families of arrows of **FinLim**. We expect to make this explicit in a later work.

M. Makkai [1993a], [1993b] has produced an approach to explicating the logic of sketches that is quite different from that presented here. Both are attempts at codifying the process of diagram-manipulation used to prove results valid in particular structured categories. In both cases, structured categories or doctrines form the semantic universes with a pre-existing notion of validity. Both approaches are motivated by the desire to formulate a syntactic notion of deducibility. After that is said, the two approaches are very different and a detailed investigation of the relationship between them is desirable.

The following points are however worth mention.

- Both approaches require a generalization of Ehresmann sketches: Sketches as in Definition 5.3.1 here and sketch category over a category $G$ in [Makkai, 1993b].

- In Makkai's approach, the sketch axioms are different for different doctrines and serve both as axioms and as rules of inference. In contrast, here the rules of construction are the *same* for all doctrines, what distinguishes doctrines is their specification as CS-sketches. This feature is a departure from the usual practice in symbolic logic.

The rules of construction given herein take place in the doctrine of finite limits. Most of the syntax and rules of deduction of traditional logical theories are clearly expressible using context-sensitive grammars, which intuitively at least can be modeled using finite limits. (Context-free grammars can be modeled using only finite products [Wells and Barr, 1988].) However, one could imagine extensions of this doctrine; for example, the use of coproducts to allow the specification of conditional compilation or the use of doctrines involving epimorphic covering families that allow the deliberate description of ambiguous statements.

## A  Appendix: Rules of construction

These rules are given in two lists.

### A.1  Rules that construct objects and arrows

This list gives all the rules that construct objects and arrows in the category CatTh[**FinLim**, S].

$$\exists\text{OB} \quad \frac{\phantom{xxx}}{c} \quad \text{for every object } c \text{ of } S.$$

$$\exists\text{ARR} \quad \frac{\phantom{xxxxxx}}{a \xrightarrow{f} b} \quad \text{for every arrow } f : a \to b \text{ of } S.$$

$$\exists\text{COMP} \quad \frac{a \xrightarrow{f} b \xrightarrow{g} c}{a \xrightarrow{g \circ f} c} \quad \text{for every object } b \text{ and pair of arrows } f : a \to b \text{ and } g : b \to c \text{ of CatTh[\textbf{FinLim}, S]}.$$

$$\exists\text{ID} \quad \frac{c}{\mathsf{Id}[c] \;\circlearrowright c} \quad \text{for every object } c \text{ of } \mathsf{CatTh}[\mathbf{FinLim}, S].$$

$$\exists\text{LIM} \quad \frac{\delta : I \to \mathsf{CatTh}[\mathbf{FinLim}, S]}{\mathsf{LimCone}[\delta] : \mathsf{Lim}[\delta] \lhd \delta} \quad \text{for every diagram } \delta : I \to \mathsf{CatTh}[\mathbf{FinLim}, S].$$

$$\exists\text{FIA} \quad \frac{\Theta : v \lhd \delta}{\mathsf{Fillin}[\Theta, \delta] : v \to \mathsf{Lim}[\delta]} \quad \text{for every diagram } \delta \text{ and every cone } \Theta : v \lhd \delta \text{ in } \mathsf{CatTh}[\mathbf{FinLim}, S].$$

**A.1.1 Remark** The first two rules are justified by the inclusion of the sketch $S$ into $\mathsf{CatTh}[\mathbf{FinLim}, S]$. The rules ∃COMP corresponds to the arrow comp and ∃ID to unit of the sketch for categories B.2. ∃LIM and ∃FIA correspond to arrows in **FinLim** but not specifically to arrows of the sketch B.5, because an arbitrary finite limit is constructed from a combination of products and equalizers.

### A.1.2 Remark

The rules just given construct specific objects and arrows in $\mathsf{CatTh}[\mathbf{FinLim}, S]$. Rule ∃LIM, for example, constructs a specific limit cone called $\mathsf{LimCone}[\delta]$, thus providing canonical limits for $\mathsf{CatTh}[\mathbf{FinLim}, S]$. It is true that there are other limit cones in general for a given diagram $\delta$, but $\mathsf{LimCone}[\delta]$ is a specific one.

Of course, in many cases, the entity constructed is the unique entity satisfying some property. For example, the arrow $g \circ f$ constructed by ∃COMP is (by definition of commutative diagram) the only one making the bottom diagram in COMPDIAG commute in $\mathsf{CatTh}[\mathbf{FinLim}, S]$. The arrow constructed by ∃ID is (by an easy theorem of category theory) the only one making the bottom diagrams in IDL and IDR commute. The arrow constructed by ∃FIA is (because of !FIA) the only one making the bottom diagram in CFIA commute. In connection with the point that each rule constructs a specific arrow, these observations are red herrings: in fact, each rule constructs a specific arrow with the name given, independently of any uniqueness properties arising from any other source.

This point of view is contrary to the spirit of category theory. We follow it here because we are constructing syntax with an eye toward implementing it in a computer language. It may be possible to design a computer language in which a construction such as ∃LIM produces a pointer to an isomorphism class of structures rather than to a specific one of them, but we know of no such language. On the other hand, implementing the specific constructions defined above should be relatively straightforward using a modern object-oriented language. Note that we are not asserting that it would be straightforward to find a confluent and normalizing form of these rules for automatic theorem proving, only that there are no obvious difficulties in implementing them so that they could be applied in an *ad-hoc* manner.

## A.2 Rules that construct diagrams

The following rules produce the existence of diagrams that must commute.

$$\text{REF} \quad \frac{a \xrightarrow{f} b}{a \underset{f}{\overset{f}{\rightrightarrows}} b} \quad \text{for every arrow } f : a \to b \text{ of } \mathsf{CatTh}[\mathbf{FinLim}, S]$$

$$\text{TRANS} \quad \frac{a \underset{g}{\overset{f}{\rightrightarrows}} b \quad a \underset{h}{\overset{g}{\rightrightarrows}} b}{a \underset{h}{\overset{f}{\rightrightarrows}} b} \quad \begin{array}{l} \text{for all objects } a \text{ and } b \text{ and} \\ \text{all arrows } f, g, h : a \to b \text{ of} \\ \mathsf{CatTh}[\mathbf{FinLim}, S] \end{array}$$

$$\exists\text{DIAG} \quad \frac{}{\mathsf{UnivMod}[\mathbf{FinLim}, S] \circ \delta : I \to \mathsf{CatTh}[\mathbf{FinLim}, S]}$$

for every diagram $\delta$ in the set $D_S$ of diagrams of $S$.

| | | |
|---|---|---|
| COMPDIAG | $$\dfrac{\begin{array}{c}a \xrightarrow{f} b \\ \downarrow g \\ c\end{array}}{\begin{array}{c}a \xrightarrow{f} b \\ {}_{g\circ f}\searrow \;\downarrow g \\ c\end{array}}$$ | for every pair of arrows $f\colon a \to b$ and $g : b \to c$ of $\mathsf{CatTh}[\mathbf{FinLim}, S]$. |
| IDL | $$\dfrac{\begin{array}{c}\mathsf{Id}[c]\circlearrowleft c\end{array}\quad g : b \to c}{\begin{array}{c}b \xrightarrow{g} c \\ {}_{g}\searrow \;\downarrow \mathsf{Id}[c] \\ c\end{array}}$$ | for every object $c$ and every arrow $g : b \to c$ of $\mathsf{CatTh}[\mathbf{FinLim}, S]$. |
| IDR | $$\dfrac{\begin{array}{c}\mathsf{Id}[c]\circlearrowleft c\end{array}\quad f : c \to d}{\begin{array}{c}c \xrightarrow{\mathsf{Id}[c]} c \\ {}_{f}\searrow \;\downarrow f \\ d\end{array}}$$ | for every object $c$ and every arrow $f : c \to d$ of $\mathsf{CatTh}[\mathbf{FinLim}, S]$. |
| ASSOC | $$\dfrac{a \xrightarrow{f} b \xrightarrow{g} c \xrightarrow{h} d}{\begin{array}{c}a \xrightarrow{f} b \\ {}_{g\circ f}\downarrow \;\swarrow g\; \downarrow h\circ g \\ c \xrightarrow{h} d\end{array}}$$ | for all arrows $f : a \to b$, $g : b \to c$ and $h : c \to d$ of $\mathsf{CatTh}[\mathbf{FinLim}, S]$. |

CFIA

$$\frac{i \in \mathsf{Nodes}[I] \qquad \Theta : v \triangleleft \delta}{}$$

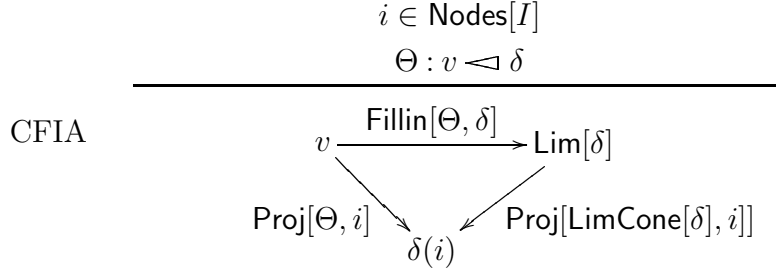

for every diagram $\delta\colon I \to \mathsf{CatTh}[\mathbf{FinLim}, S]$, every node $i$ of $I$, and every cone $\Theta$ with base diagram $\delta$.

!FIA

$$\delta : I \to \mathsf{CatTh}[\mathbf{FinLim}, S]$$
$$\Theta : v \triangleleft \delta$$
$$h : v \to \mathsf{Lim}[\delta]$$
$$k : v \to \mathsf{Lim}[\delta]$$

and each of the following diagrams for each node $i$ of $I$:

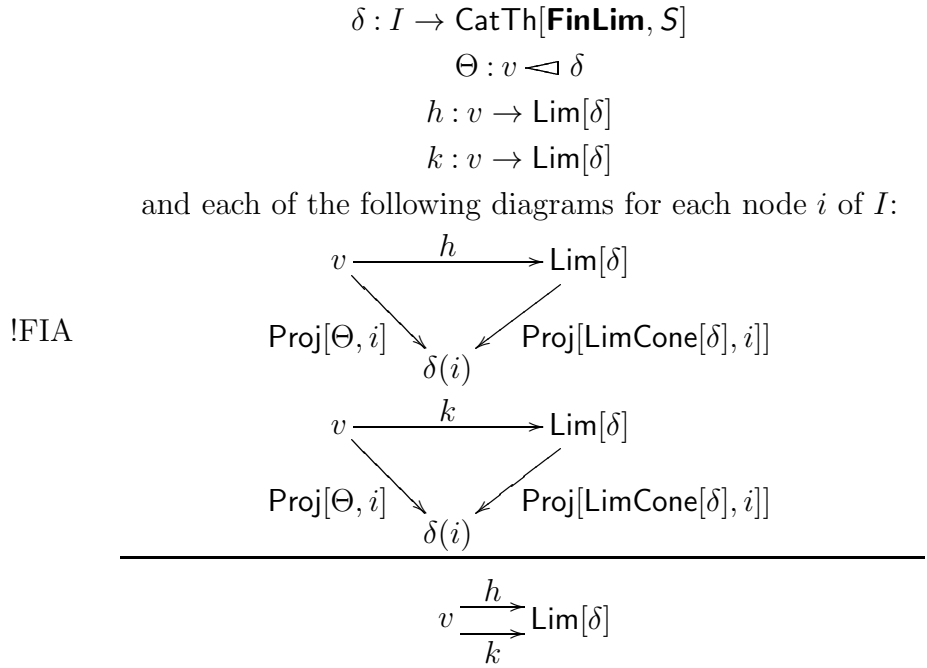

for every diagram $\delta\colon I \to \mathsf{CatTh}[\mathbf{FinLim}, S]$, every cone $\Theta$ in $\mathsf{CatTh}[\mathbf{FinLim}, S]$ with base diagram $\delta$, and every pair of arrows $h, k : v \to \mathsf{Lim}[\delta]$.

**A.2.1 Remark** Note that we do not need a rule of the form

$$\text{SYM} \quad \frac{a \underset{g}{\overset{f}{\rightrightarrows}} b}{a \underset{f}{\overset{g}{\rightrightarrows}} b}$$

since the two diagrams exhibited are actually the same diagram (see 2.3).

# B  Appendix: Sketches for constructor spaces

Here we present constructor space sketches for certain types of categories. In each case the models are categories of the sort described and the morphisms of models are functors that preserve the structure on the nose. It is an old result that such categories can be sketched. See [Burroni, 1970a], [Burroni, 1970b], [McDonald and Stone, 1984], and [Coppey and Lair, 1988], for example.

The embedding $\eta$ of Section 5.1.3 will in each case be inclusion.

## B.1  Notation

We denote the $i$th projection in a product diagram of the form

$$\begin{array}{c} \mathsf{ob}^A \times \mathsf{ob}^B \\ {}^{p_1}\swarrow \qquad \searrow^{p_2} \\ \mathsf{ob}^A \qquad\qquad \mathsf{ob}^B \end{array}$$

as $p_i$, or $p_i^{A\times B}$ if the source or target is not shown. We use a similar device for the product of three copies of $\mathsf{ob}$.

## B.2  The sketch *Cat* for categories

This version of the sketch for categories is based on [Barr and Wells, 1990]. Another version is given in [Coppey and Lair, 1988], page 64. The first versions were done by Ehresmann [1966], [1968a] and [1968b].

### B.2.1 The graph of *Cat*

The graph of the sketch for categories contains nodes as follows.

1. **1**, the formal terminal object.
2. ob, the formal set of objects.
3. ar, the formal set of arrows.
4. $\mathsf{ar}_2$, the formal set of composable pairs of arrows.
5. $\mathsf{ar}_3$, the formal set of composable triples of arrows.

The arrows for the sketch for categories are

1. unit : ob $\to$ ar that formally picks out the identity arrow of an object.
2. source, target : ar $\to$ ob that formally pick out the source and target of an arrow.
3. comp : $\mathsf{ar}_2 \to$ ar that picks out the composite of a composable pair.
4. lfac, rfac: $\mathsf{ar}_2 \to$ ar that pick out the left and right factors in a composable pair.
5. lfac, mfac, rfac : $\mathsf{ar}_3 \to$ ar that pick out the left, middle and right factors in a composable triple of arrows.
6. lass, rass : $\mathsf{ar}_3 \to \mathsf{ar}_2$: lass formally takes $\langle h, g, f \rangle$ to $\langle h \circ g, f \rangle$ and rass takes it to $\langle h, g \circ f \rangle$.
7. lunit, runit : ar $\to \mathsf{ar}_2$: lunit takes an arrow $f : A \to B$ to $\langle \mathsf{Id}[B], f \rangle$ and runit takes it to $\langle f, \mathsf{Id}[A] \rangle$.
8. Arrows id : $x \to x$ as needed.

Observe that id, lfac and rfac, like $p_1$ and $p_2$, are overloaded. We will observe the same care with these arrows as with $p_1$ and $p_2$ as mentioned in Section B.1.

### B.2.2 Cones of *Cat*

$\text{ar}_2$ and $\text{ar}_3$ are defined by these cones:

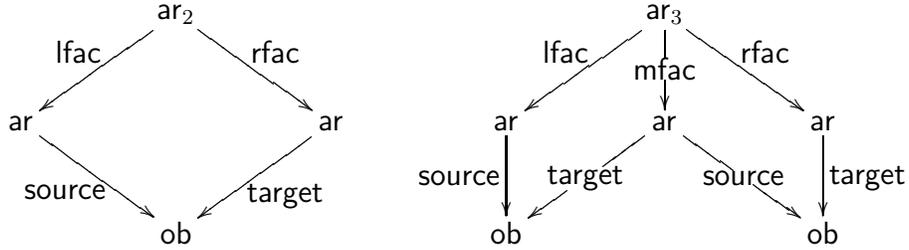

### B.2.3 Diagrams of *Cat*

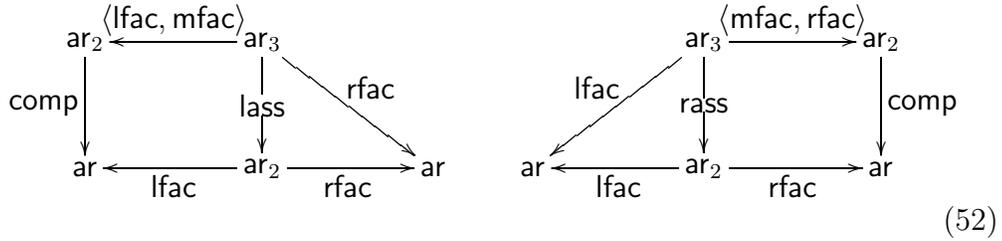

(52)

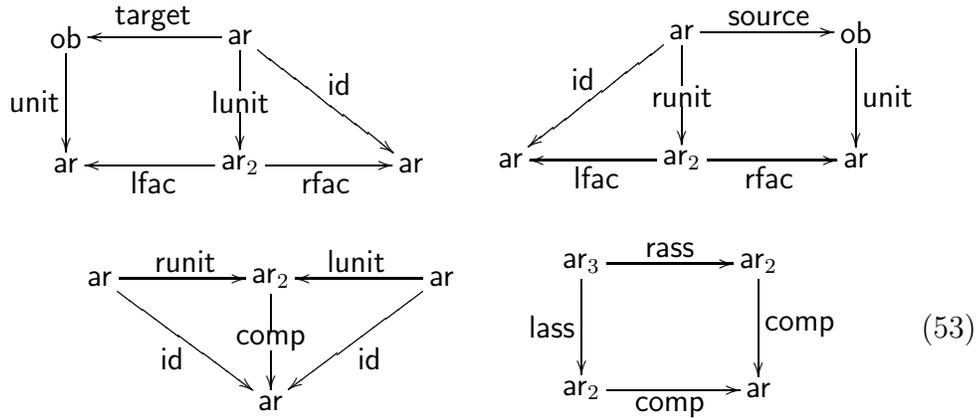

(53)

### B.3 The sketch for the constructor space **FinProd**

To get the sketch for categories with finite products, we must add the following nodes and arrows to the sketch for categories:

Nodes:

1. ta, the formal set of terminal arrows.

2. cone, the formal set of cones of the form

$$\begin{array}{c} & V & \\ {}^{p_1}\swarrow & & \searrow^{p_2} \\ A & & B \end{array} \qquad (54)$$

3. fid, the formal set of fill-in diagrams ("sawhorses") of the form

$$\begin{array}{ccc} V & \xrightarrow{h} & L \\ \downarrow & \times & \downarrow \\ A & & B \end{array} \qquad (55)$$

where $h$ commutes with the cone projections.

Arrows:

1. ter : $1 \to$ ob, that formally picks out a particular terminal object.

2. ! : ob $\to$ ta, that picks out the arrow from an object to the terminal object.

3. inc : ta $\to$ ar, the formal inclusion of the set of terminal arrows into the set of arrows.

4. prod : ob $\times$ ob $\to$ cone, that picks out the product cone over a pair of objects.

5. soco : fid $\to$ cone, that picks out the source cone of a fill-in arrow.

6. taco : fid $\to$ cone, that picks out the target cone of a fill-in arrow.

7. ufid : cone $\to$ fid, that takes a cone to the unique fill-in diagram that has the cone as source cone.

8. fia : fid $\to$ ar that formally picks out the fill-in arrow in a fill-in diagram.

### B.3.1 Cones for **FinProd**

**FinProd** has four cones in addition to those of the sketch for categories. One is the cone with vertex 1 over the empty diagram. The one below says that ta is the formal set of arrows to the terminal object:

$$
\begin{array}{c}
\text{ta} \\
\text{inc} \swarrow \quad \searrow \\
\text{ar} \xrightarrow{\text{target}} \text{ob} \\
\searrow \quad \nearrow \text{ter} \\
1
\end{array}
\tag{56}
$$

Note that in giving this cone, we are not only saying that ta is the limit of the diagram

$$
\text{ar} \xrightarrow{\text{target}} \text{ob} \nwarrow_{\text{ter}} 1
\tag{57}
$$

but also that inc is one of the projection arrows. (Indeed, this is the only projection arrow that matters, since the other two are induced.)

The following cone makes cone the formal object of cones to a discrete diagrams consisting of a pair of objects.

$$
\begin{array}{c}
\text{cone} \\
\text{lproj} \swarrow \quad \searrow \text{rproj} \\
\text{ar} \xrightarrow{\text{source}} \text{ob} \xleftarrow{\text{source}} \text{ar}
\end{array}
\tag{58}
$$

Finally, there must be a cone with vertex fid over Diagram 59 below, which is annotated to refer to Diagram (55), in which $\Gamma$ is the cone with vertex $V$, $\Lambda$ is the cone with vertex $L$ and is a limit cone, and $h$ is the fill-in arrow. In this case, the projection arrows of the cone are not shown.

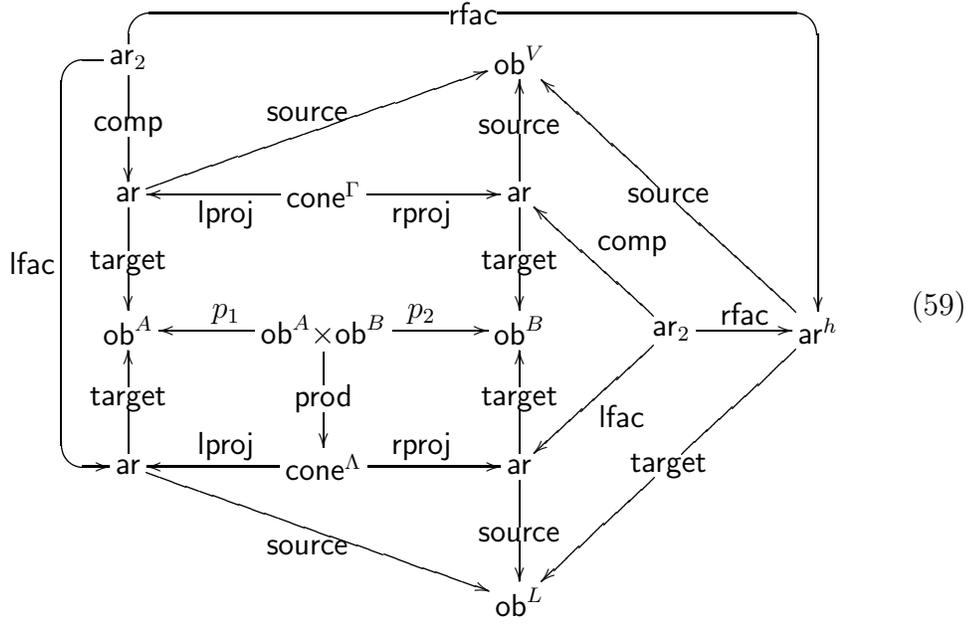

(59)

In addition, we require:

1) The projection to $\mathsf{cone}^\Gamma$ must be $\mathsf{soco}$.
2) The projection to $\mathsf{cone}^\Lambda$ must be $\mathsf{taco}$.
3) The projection to $\mathsf{ar}^h$ must be $\mathsf{fia}$.

### B.3.2  Diagrams for **FinProd**

The following two diagrams make the arrow to the terminal object have the correct source and target.

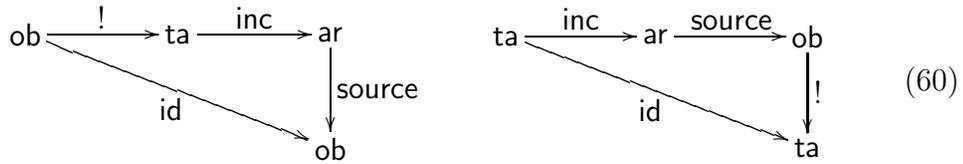

(60)

The diagram below makes the fill-in arrow to a product unique.

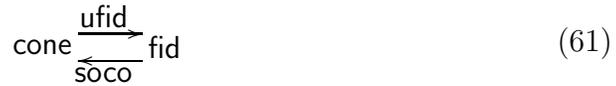

(61)

The diagram below forces the product cone projections to have the correct targets.

$$
\begin{array}{c}
\text{ob} \times \text{ob} \\
{}_{p_1}\swarrow \quad \downarrow \text{prod} \quad \searrow {}^{p_2} \\
\text{ob} \xleftarrow{\text{target}} \text{ar} \xleftarrow{\text{lproj}} \text{cone} \xrightarrow{\text{rproj}} \text{ar} \xrightarrow{\text{target}} \text{ob}
\end{array}
\tag{62}
$$

## B.4  Modules

As we proceed to sketch more complicated constructions, we will need to use some device to communicate the nature of the necessary diagrams, which become too large to comprehend easily. Here we introduce the first of several **modules**: diagrams that occur frequently as subdiagrams because they are needed to force the value of a node in a model to contain certain types of constructions.

### B.4.1  The module for the product of objects

Every occurrence of ob that is annotated $M \times N$ must be part of a subdiagram of the following form:

$$
\begin{array}{ccc}
\text{ob}^M & \xleftarrow{\text{target}} & \text{ar}^{p_1^{M \times N}} \\
{}^{p_1}\uparrow & & \uparrow \text{lproj} \quad \searrow {}^{\text{source}} \\
\text{ob}^M \times \text{ob}^N & \xrightarrow{\text{prod}} & \text{cone} \quad\quad \text{ob}^{M \times N} \\
{}^{p_2}\downarrow & & \downarrow \text{rproj} \quad \nearrow {}^{\text{source}} \\
\text{ob}^N & \xleftarrow{\text{target}} & \text{ar}^{p_2^{M \times N}}
\end{array}
\tag{63}
$$

Henceforth, an occurrence of ob annotated $A \times B$ (for example) will be taken to *imply* the existence of a subdiagram of the form of Diagram (63) with $M$ replaced with $A$ and $N$ replaced with $B$. The subdiagram will not necessarily be shown. If this is part of a diagram $\delta$, the diagram can be reconstructed by taking the union of the shape graph of the module (63) and the shape graph of the part of $\delta$ that is shown on the page, and defining the diagram based on the resulting graph as the pushout of the diagram shown and the module. This is illustrated in Diagrams (65) and (66) in the next section.

### B.4.2 The module for the product of arrows

In the commutative diagram

$$\begin{array}{ccc} K & \xleftarrow{p_1} & K \times N \\ & & \\ u \downarrow & \downarrow & \searrow^{p_2} \\ & & N \\ & & \nearrow_{p_2} \\ M & \xleftarrow{p_2} & M \times N \end{array} \qquad (64)$$

the unlabeled arrow is necessarily $u \times \mathsf{Id}[N] : K \times N \to M \times N$. Such a diagram must be an element in a model of the value of Diagram (65) below, which is therefore a module for the product of an arrow and an identity arrow. In this diagram, $\phi := \langle p_2^{M \times N}, u \times \mathsf{Id}[N] \rangle$.

$$(65)$$

More precisely, let $x$ be an element of $\mathfrak{M}(\mathsf{Lim}[D(65)])$, for some category $\mathfrak{M}$ with finite limits. Then if $\mathsf{Proj}[\mathsf{Lim}[D(65)], h](x) = h$, then

$$\mathsf{Proj}[\mathsf{Lim}[D(65)], h\times \mathsf{Id}[A]](x) = h\times \mathsf{Id}[A]$$

as suggested by the notation. This will be used in Section B.6 below.

Diagram (65) contains two copies of Diagram (63), the module for the product of two objects. The copy at the bottom is precisely Diagram (63), and the copy at the top is Diagram (63) with $M$ replaced with $K$. In the sequel, a diagram such as Diagram (65) will be drawn without the modules, as shown below.

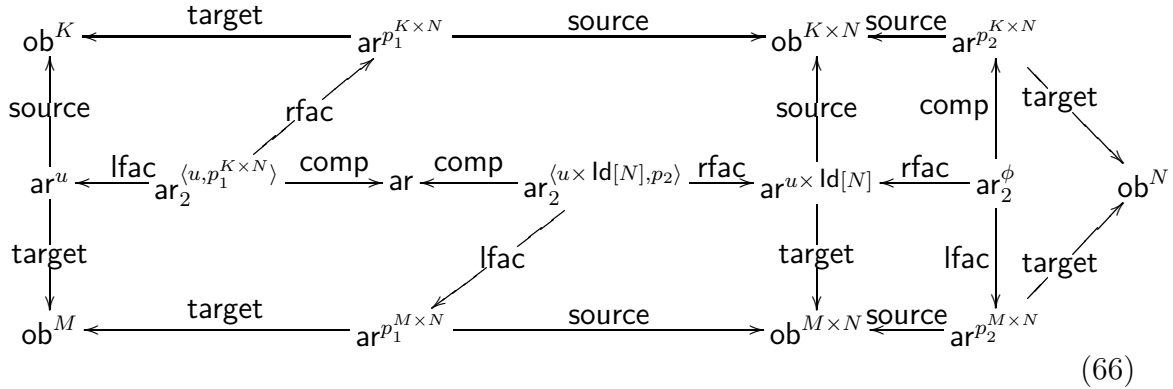
(66)

Diagram (65) may be mechanically reconstructed from Diagram (66) and the annotations that include the symbols $M \times N$ and $K \times N$ (three of each). The shape diagram of Diagram (65) is the pushout of the shape diagram of Diagram (66) and the shape diagrams of the modules Diagram (63) and Diagram (63) with $M \leftarrow K$. Each of the latter two have four annotated nodes and six annotated arrows in common with Diagram (66), and the values of any two of the three smaller diagrams at a given common node or arrow is of course the same, so that Diagram (65) is the union of Diagram (66) and the two modules.

## B.5 The sketch for the constructor space **FinLim**

We sketch the constructor space **FinLim** by adding data to the sketch for **FinProd** that ensure that a **FinLim**-category has equalizers of pairs of arrows. The sketch has the following nodes:

1. ppair is the formal set of parallel pairs of the form

$$A \xrightarrow[g]{f} B \tag{67}$$

2. econe is the formal set of diagrams

$$E \xrightarrow{u} A \xrightarrow[g]{f} B \tag{68}$$

in which $f \circ u = g \circ u$. Of course, a cone to Diagram (67) also has a projection to $B$, but that is forced and need not be included in the data for the cone.

3. efid is the set of fill-in diagrams

$$\begin{array}{c} X \\ {}^{v}\swarrow \phantom{x} \downarrow^{u} \\ E \xrightarrow{e} A \xrightarrow[g]{f} B \end{array} \tag{69}$$

in which $f \circ e = g \circ e$ and $u = e \circ v$.

The arrows of the sketch include:

1. equ : ppair $\to$ econe, that formally picks out the equalizer of the parallel pair.

2. top, bot : ppair $\to$ ar, that pick out $f$ and $g$ in Diagram (67).

3. etop, ebot : econe $\to$ ar, that pick out $f$ and $g$ in Diagram (68).

4. esoco, etaco : efid $\to$ econe that pick out the source and target cones of the fill-in arrow.

5. eufid : econe $\to$ efid that takes a diagram of the form of Diagram (68) to the unique fill-in diagram that has this diagram as source cone.

6. efia : efid $\to$ ar that picks out the fill-in arrow in a fill-in diagram.

### B.5.1 Cones for **FinLim**

ppair is the limit of the diagram

$$
\begin{array}{ccc}
\mathsf{ar}^f & \xrightarrow{\mathsf{target}} & \mathsf{ob}^B \\
{\scriptstyle\mathsf{source}}\downarrow & & \uparrow{\scriptstyle\mathsf{target}} \\
\mathsf{ob}^A & \xleftarrow{\mathsf{source}} & \mathsf{ar}^g
\end{array}
\tag{70}
$$

The following projections from ppair have names: $\mathsf{soob} : \mathsf{ppair} \to \mathsf{ob}^A$, $\mathsf{top} : \mathsf{ppair} \to \mathsf{ar}^f$, and $\mathsf{bot} : \mathsf{ppair} \to \mathsf{ar}^g$.

econe is the limit of

$$
\text{(diagram)} \tag{71}
$$

Two projections have names: $\mathsf{etop} : \mathsf{econe} \to \mathsf{ar}^f$ and $\mathsf{ebot} : \mathsf{econe} \to \mathsf{ar}^g$.

efid is the limit of the pushout of Diagram (71) and the following diagram. Note that the common part of the two diagrams is

$$
\mathsf{ar}^e \xrightarrow{t} \mathsf{ob}^A
$$

We could have presented Diagram (59) as a pushout in much the same way (the common part would describe the arrow $h : V \to L$). We have deliberately varied the way we present the data in this article because we are not sure ourselves which approach communicates best.

$$
\text{(diagram)} \tag{72}
$$

The named projections are esoco: efid → $\text{ar}^u$, etaco: efid → $\text{ar}^e$ and efia: efid → $\text{ar}^v$.

### B.5.2   Diagrams for **FinLim**

The following diagram makes the fill-in arrow unique.

$$\text{econe} \underset{\text{esoco}}{\overset{\text{eufid}}{\rightleftarrows}} \text{efid} \tag{73}$$

These two diagrams ensure that the equalizer cone be a cone to the correct diagram.

$$\begin{array}{c}\text{ppair} \xrightarrow{\text{equ}} \text{econe} \\ \scriptstyle{\text{top}} \searrow \downarrow \scriptstyle{\text{etop}} \\ \text{ar}\end{array} \qquad \begin{array}{c}\text{ppair} \xrightarrow{\text{equ}} \text{econe} \\ \scriptstyle{\text{bot}} \searrow \downarrow \scriptstyle{\text{ebot}} \\ \text{ar}\end{array} \tag{74}$$

## B.6   The sketch for the constructor space **CCC**

### B.6.1   Definition
A **Cartesian closed category** is a category $\mathscr{C}$ with the following structure:

CCC.1  $\mathscr{C}$ has binary products.

CCC.2  For each pair of objects $A$ and $B$ of $\mathscr{C}$, there is an object $B^A$ and an arrow eval: $B^A \times A \to B$.

CCC.3  For each triple of objects $A$, $B$ and $C$ of $\mathscr{C}$, there is a map

$$\lambda : \text{Hom}(B \times A, C) \to \text{Hom}(B, C^A) \tag{75}$$

such that for every arrow $f : B \times A \to C$,

$$\begin{array}{c} B \times A \xrightarrow{\lambda f \times \text{Id}[A]} C^A \times A \\ \scriptstyle{f} \searrow \qquad \downarrow \scriptstyle{\text{eval}} \\ C \end{array} \tag{76}$$

commutes.

CCC.4  For any arrow $g : B \to C^A$, $\lambda(\text{eval} \circ (g \times \text{Id}[A])) = g$.

Using this definition, the sketch for the constructor space for Cartesian closed categories may be built on the sketch for **FinProd** by adding the following nodes and arrows.

The nodes are:

1. twovf, the formal set of "functions of two variables", that is, arrows of the form $B \times A \to C$.

2. curry, the formal set of "curried functions" $B \to C^A$.

The sketch for **CCC** has arrows

1. $\mathsf{fs} : \mathsf{ob}^B \times \mathsf{ob}^A \to \mathsf{ob}^{B^A}$ that picks out the function space $B^A$ of two objects $B$ and $A$.

2. $\mathsf{ev} : \mathsf{ob}^B \times \mathsf{ob}^A \to \mathsf{ar}$ that picks out the arrow $\mathsf{eval} : B^A \times A \to B$.

3. $\mathsf{lam} : \mathsf{ar} \to \mathsf{ar}$, the formal version of the mapping $\lambda$ of Diagram (75).

4. $\mathsf{tsource} : \mathsf{twovf} \to \mathsf{ob}^{B \times A}$, that picks out the source of a function $f : A \times B \to C$.

5. $\mathsf{ttarget} : \mathsf{twovf} \to \mathsf{ob}^C$, that picks out the target of a function $f : A \times B \to C$.

6. $\mathsf{arrow} : \mathsf{twovf} \to \mathsf{ar}^f$, that picks out the arrow $f$ itself.

7. $\mathsf{csource} : \mathsf{curry} \to \mathsf{ob}^B$, that picks out the source of a curried function $g : B \to C^A$.

8. $\mathsf{ctarget} : \mathsf{curry} \to \mathsf{ob}^{C^A}$, that picks out the target of a curried function $f : B \to C^A$.

9. $\mathsf{arrow} : \mathsf{curry} \to \mathsf{ar}^g$, that picks out the arrow $g$ itself.

### B.6.2 Cones for **CCC**

**CCC** must have two cones

$$\begin{array}{c} \text{twovf} \\ \swarrow \;\; \downarrow \;\; \searrow \\ \text{tsource} \quad \text{arrow} \quad \text{ttarget} \\ \mathsf{ob}^{B \times A} \xleftarrow{\text{source}} \mathsf{ar}^f \xrightarrow{\text{source}} \mathsf{ob}^C \end{array} \quad (77)$$

$$\begin{array}{c}
\text{curry} \\
\swarrow \quad \downarrow \text{ctarget} \quad \text{arrow} \searrow \text{csource} \\
\mathsf{ob}^C \times \mathsf{ob}^A \xrightarrow{\text{fs}} \mathsf{ob}^{c^A} \xleftarrow{\text{target}} \mathsf{ar}^g \xrightarrow{\text{source}} \mathsf{ob}^B
\end{array} \qquad (78)$$

## B.7  The module for function spaces

Henceforth, we will assume the module

$$\mathsf{ob}^M \xleftarrow{p_1} \mathsf{ob}^M \times \mathsf{ob}^N \xrightarrow{\text{fs}} \mathsf{ob}^{M^N} \\ \downarrow{p_2} \\ \mathsf{ob}^N \qquad (79)$$

is attached whenever an occurrence of $\mathsf{ob}$ is annotated by $M^N$. Note that this occurred in Diagram (78).

### B.7.1  Diagrams for CCC

Diagram (80) below forces eval to have the correct domain and codomain.

$$\mathsf{ob}^{B^A \times A} \xleftarrow{\text{source}} \mathsf{ar}^{\text{eval}} \xrightarrow{\text{target}} \mathsf{ob}^B \qquad (80)$$

Expanding this diagram using the required modules is a two stage process, giving

$$\begin{array}{c} \text{(diagram)} \end{array} \qquad (81)$$

Diagram (82) below forces $\lambda f$ to have the correct domain and codomain.

$$\begin{array}{ccc}
& \mathrm{ob}^C & \mathrm{ob}^{C^A} \\
& \uparrow\text{ttarget} & \uparrow\text{csource} \\
\mathrm{ob}^{B\times A} \xleftarrow{\text{tsource}} \text{twovf}^f \xrightarrow{\text{lam}} & \text{curry}^{\lambda f} \\
& & \downarrow\text{ctarget} \\
& & \mathrm{ob}^B
\end{array} \qquad (82)$$

Diagram (83) below forces Diagram (76) to commute.

$$\begin{array}{c}
\mathrm{ob}^{B\times A} \xleftarrow{\text{source}} \mathrm{ar}^{\lambda f \times \mathrm{Id}[A]} \xrightarrow{\text{target}} \mathrm{ob}^{C^A \times A} \\
\uparrow \text{source} \quad \uparrow\text{rfac} \quad \uparrow\text{source} \\
\quad \mathrm{ar}_2 \\
\swarrow\text{comp} \qquad \searrow\text{lfac} \\
\mathrm{ar}^f \xrightarrow{\text{target}} \mathrm{ob}^C \xleftarrow{\text{target}} \mathrm{ar}^{\text{eval}}
\end{array} \qquad (83)$$

Diagram (84) below ensures that requirement CCC–4 holds.

$$\begin{array}{c}
\mathrm{ob}^{C^A\times A} \xleftarrow{\text{source}} \mathrm{ar}^{\text{eval}} \xrightarrow{\text{target}} \mathrm{ob}^C \qquad \mathrm{ob}^B \\
\uparrow\text{target} \quad \uparrow\text{lfac} \qquad \qquad \uparrow\text{source} \\
\mathrm{ar}^{g\times\mathrm{Id}[A]} \xleftarrow{\text{rfac}} \mathrm{ar}_2^{\langle\text{eval},g\times\mathrm{Id}[A]\rangle} \quad \text{target} \quad \mathrm{ar}^g \xrightarrow{\text{target}} \mathrm{ob}^{C^A} \\
\downarrow\text{source} \qquad \searrow\text{comp} \quad \uparrow\text{ttarget} \quad \uparrow\text{arrow} \\
\mathrm{ob}^{B\times A} \xleftarrow{\text{source}} \mathrm{ar}^{\text{eval}\circ(g\times\mathrm{Id}[A])} \qquad \text{curry} \\
\qquad \text{tsource} \qquad \text{arrow} \quad \uparrow\text{lam} \\
\qquad\qquad\qquad \text{twovf}
\end{array} \qquad (84)$$

### B.7.2 Invertibility of $\lambda$

It follows from CCC-4 that if $\mathscr{C}$ is any Cartesian closed category corresponding to a model $\mathfrak{C}$, then $\mathfrak{C}(\lambda)$ is a bijection. The Completeness Theorems 6.2.1 and 6.2.3 then imply that there is an arrow $\lambda^{-1}$ in **CCC** that, as its name suggests, is a formal inverse to $\lambda$.

# References

At the end of each entry, the pages on which that entry is cited are listed in parentheses.

Atish Bagchi
226 West Rittenhouse Square
#702
Philadelphia, PA 19103
atish@math.upenn.edu

Charles Wells
Department of Mathematics
Case Western Reserve University
University Circle
Cleveland, OH 44106-7058, USA
cfw2@po.cwru.edu